\documentclass[12pt,a4paper]{article}

\usepackage{a4wide,amssymb,amsmath,amsthm,xspace,epsfig, amsfonts,amscd}
\usepackage{tikz-cd}

\usepackage[left=1.5cm, right=1.5cm, top=1.5cm]{geometry}

\begin{document}

\newcommand{\N}{\mathbb{N}}
\newcommand{\R}{\mathbb{R}}
\newcommand{\Z}{\mathbb{Z}}
\newcommand{\Q}{\mathbb{Q}}
\newcommand{\C}{\mathbb{C}}
\newcommand{\PP}{\mathbb{P}}

\newcommand{\VV}{\mathsf{V}}
\newcommand{\HH}{\mathsf{H}}

\newcommand{\LL}{\Bbb L}
\newcommand{\OO}{\mathcal{O}}

\newcommand{\BO}{\text{Bo}}

\newcommand{\esp}{\vskip .3cm \noindent}
\mathchardef\flat="115B

\newcommand{\lev}{\text{\rm Lev}}

\newcommand{\co}{\colon}

\def\ut#1{$\underline{\text{#1}}$}
\def\CC#1{${\cal C}^{#1}$}
\def\h#1{\hat #1}
\def\t#1{\tilde #1}
\def\wt#1{\widetilde{#1}}
\def\wh#1{\widehat{#1}}
\def\wb#1{\overline{#1}}

\def\restrict#1{\bigr|_{#1}}

\def\hu#1#2{\mathsf{U}_{fin}\bigl({#1},{#2}\bigr)}
\def\ch#1#2{\left(\begin{array}{c}#1 \\ #2 \end{array}\right)}

\newtheorem{lemma}{Lemma}[section]

\newtheorem{thm}[lemma]{Theorem}
\newtheorem*{thm*}{Theorem}
\newtheorem*{lemma?}{Lemma ??}

\newtheorem{defi}[lemma]{Definition}
\newtheorem{conj}[lemma]{Conjecture}
\newtheorem{cor}[lemma]{Corollary}
\newtheorem{prop}[lemma]{Proposition}
\newtheorem{prob}[lemma]{Problem}
\newtheorem{qu}[lemma]{Question}
\newtheorem{q}[lemma]{Question}
\newtheorem{sq}{SELF QUESTION}

\newtheorem{examples}[lemma]{Examples}
\newtheorem{example}[lemma]{Example}

\theoremstyle{remark}
\newtheorem{claim}{Claim}[lemma]
\renewcommand{\theclaim}{\arabic{section}.\arabic{lemma}.\arabic{claim}}

\newtheorem*{rem}{Remark}
\newtheorem{rem_numbered}[lemma]{Remark}

\title{Relative (functionally) Type I spaces and narrow subspaces}
\date{August 22, 2022}
\author{Mathieu Baillif}
\maketitle

\abstract{\footnotesize
An open chain cover $\mathcal{U}=\{U_\alpha\,:\,\alpha\in\kappa\}$ ($\kappa$ a cardinal) of a space $X$
is a systematic cover
if $\wb{U_\alpha}\subset U_\beta$ when $\alpha<\beta$, and
$X$ is Type I if $\kappa=\omega_1$
and each $\wb{U_\alpha}$ is Lindel\"of.
A closed subspace $D\subset X$ is narrow in $X$ if for each systematic cover $\{V_\alpha\,:\,\alpha\in\omega_1\}$ of $X$,
either there is $\alpha$ such that $D\subset V_\alpha$, or $\wb{V_\alpha}\cap D$ is Lindel\"of for each $\alpha$.
Taking systematic covers given by $s^{-1}([0,\alpha))$ for a continuous $s\co X\to\LL_{\ge 0}$ (where $\LL_{\ge 0}$ is the longray)
defines functionally Type I spaces and functionally narrow subspaces. For instance,
$\LL_{\ge 0}$ and $\omega_1$ are narrow in themselves and any other space. \\
We investigate these properties and relative versions, as well as their relationship, and show in particular the following.
There are functionally Hausdorff Type I spaces which are not functionally Type I
while regular Type I spaces are functionally Type I.
We exhibit examples of spaces which are narrow in some but not in other spaces.
There are subspaces of a Tychonoff space $Y$ that are functionally narrow but not narrow in $Y$, 
while both notions agree if $Y$ is normal.
Under {\bf PFA} and using classical results, 
any $\omega_1$-compact locally compact countably tight Type I space contains 
a non-Lindel\"of subspace narrow in it (a copy of $\omega_1$, actually), while a Suslin tree does not.
There are spaces with subspaces narrow in them that are essentially discrete.
Finally, we investigate natural partial orders on (functionally) narrow subspaces and 
when these orders are $\omega$- or $\omega_1$-closed.}

%\tableofcontents

\section{Introduction \& definitions}\label{sec:def}

This paper is about the notions of Type I, functionally Type I, narrow and functionally narrow spaces
(and relative versions), 
and their relations.
Type I spaces were introduced by P. Nyikos \cite[2.10]{Nyikos:1984} in his study of non-metrizable manifolds.
Functionally narrow spaces were introduced by the author in the preprint \cite{mesziguessurf} 
under another
name (they were called {\em directions}, we believe that the new terminology fits better).
A good portion of
the contents of the present paper comes from another preprint \cite{mesziguesDirections} 
which is almost ten years old and has gathered a grand total of zero citations (well, now: one), thus exhibiting the intensity
of the community's interest on the subject. (Containing some unspotted inaccuracies might have been a hindrance.)
But since we keep going back to these 
questions again and again despite the passage of the years, we could not help but gather our results (old and new) in
a more polished and publishable form. Our ambition being limitless, we 
are not afraid to state publicly that we aim for a strictly positive number of readers\footnote{The boldness
of this statement is alas diminished by the fact that a failure would be without witness. Referees do not count.}.

\subsection{Informalities}

Let us try to present informally the main ideas of this paper, precise definitions can be found in the next subsection,
where the reader is encouraged to jump directly shoud they have an aversion to vagueness and shaky analogies.\\
A space is Type I (in itself) iff it is an increasing union of length $\omega_1$ of
closed Lindel\"of subspaces whose interiors contain the preceding members of the union.
For 
functionally Type I spaces, the increasing union is given by preimages of initial segments of the longray $\LL_{\ge 0}$
for a function called a slicer. 
A functionally Type I space is Type I.
If one sees spaces as human beings, a (functionally) Type I space grows slowly out of childhood (i.e. Lindel\"ofness)
and thus tends to behave more reasonably than those who 
see one aspect of their life jump at once into adulthood (going from being smooth-chinned to growing a beard overnight,
for instance), having a Lindel\"of
subspace with non-Lindel\"of closure.\footnote{We are aware that various aspects of this analogy are open to debate.}
\\
Given some ambient space $X$, a closed subspace $D$ is functionally narrow in $X$ 
if, given $f\co X\to\LL_{\ge 0}$, whenever a closed non-Lindel\"of part of $D$ is sent by $f$ into a
strict initial (Lindel\"of) segment
then all of $D$ is sent to a strict initial segment. 
Narrow subspaces are defined by a similar property, but with increasing unions of closed sets 
(whose interior contain the preceding members) 
of length $\omega_1$ of $X$ instead of a function.
A subspace that is narrow is functionally narrow.
Going back to our analogy, a narrow subspace cannot be at the
same time scolded (having a grown-up aspect treated as that of a child) 
and educated gradually into adulthood.\footnote{Well, we should maybe 
forget completely about this analogy.}
Being (functionally) narrow is very sensitive to the ambient space,
as a copy of a space can be narrow in one and non-narrow in another space. 
Prototypical examples of spaces that are narrow (in themselves and any other space) are $\omega_1$ and the longray.
\\
These concepts inspired us questions that felt natural, we try to answer some of them in the present paper.
Here is a short list.\\
-- When does being Type I imply being functionally Type I (both in a fixed ambient space)~?
   We will see that it is the case if the ambient space is normal or regular and Type I (in itself). 
   Section \ref{sec:gen} contains the proof and examples of Hausdorff spaces for which this does not hold. \\
-- When does being functionally narrow imply being narrow (both in a fixed ambient space)~?
   This is again the case if the ambient space is normal, but there are Tychonoff counterexamples.
   Section \ref{sec:narrow} is devoted to this topic.\\
-- Given a Type I space, one may ask if it contains a subspace (functionally) narrow in it.
   Counter-examples are very easy to find (for instance, the discrete space of cardinality $\omega_1$),
   but it gets trickier if one imposes global and local conditions such as countable compactness and first countability,
   in which case some results depend on the model of set theory.
   This is the subject of Section \ref{sec:without}.\\
-- Since a prototypical example of a space that is non-narrow (in itself) is the discrete space of cardinality $\omega_1$,
   it sounds fun (out of contrariness, obviously) to try to find spaces whose narrow subspaces 
   are essentially discrete. This is done in Section \ref{sec:discrete}.\\
-- The definitions of (functional) narrowness 
   yield two natural partial orders on closed non-Lindel\"of subspaces of a given space $X$.
   One is defined simply as $C\preceq_f D$ iff given any $f$ from $X$ to the longray,
   then $f$ sends $C$ into an initial segment whenever $f$ sends $D$ into an initial segment. 
   For instance, in the square of the longray, the diagonal is strictly $\preceq_f$-bigger than any horizontal or vertical
   line. These orders generalize the inclusion relation. Section \ref{sec:orders} 
   is dedicated to a short study of these orders on the set of (functionally) narrow subspaces.

\subsection{Definitions}\label{subsec:def}

By {\em space} we mean a topological Hausdorff space, hence normal and regular spaces are assumed to be Hausdorff. 
Throughout the paper, every fonction is assumed to be continuous if not stated otherwise.
We follow the set-theoretic tradition of denoting the set of positive integers by $\omega$
and the first uncountable ordinal by $\omega_1$, and $\wb{U}$ denotes
the closure of the subspace $U$ in some topological space clear from the context.
We use the term {\em club}
as a shorthand for {\em closed and unbounded}.
The restriction of a function $f$ to a subset $A$ of its domain is denoted by $f\upharpoonright A$.
We use brackets $\langle\,,\,\rangle$ for ordered pairs and parenthesis $(\,,\,)$ for open intervals
in (partially) ordered sets.
Letters $\alpha,\beta,\gamma,\delta$ are reserved exclusively for ordinals, and an ordinal is always assumed to be the
set of its predecessors. When seen as topological spaces, ordinals are endowed with the order topology (unless specified).

\begin{defi}
   A cover of a space is a cover by open sets. 
   A chain cover is linearly ordered by the inclusion relation.
   A good cover of $X$ is a chain cover $\{U_\alpha\,:\,\alpha\in\kappa\}$
   where $U_\alpha\subset U_\beta$ when $\alpha < \beta$, $\kappa$ is a regular cardinal
   and $\cup_{\gamma<\alpha}U_\gamma = U_\alpha$ when $\alpha$ is a limit ordinal.
   A systematic cover is a chain cover $\{U_\alpha\,:\,\alpha\in\kappa\}$ such that
   $\wb{U_\alpha}\subset U_\beta$  when $\alpha < \beta$.
\end{defi}

We write down the next lemma for the record since we use it throughout the paper
(most of the time implicitly). Its proof is immediate.
\begin{lemma}
   Any chain cover has a subcover which either contains only one element (the whole space)
   or is
   indexed by a regular cardinal, in the latter case it 
   can become a good chain cover by adding relevant unions at limit ordinal (when needed).
\end{lemma}

\begin{defi}
   The longray $\LL_{\ge 0}$ is $\omega_1\times[0,1)$ with lexicographic order topology.
   A subset of $\LL_{\ge 0}$ or $\omega_1$ is bounded iff it is contained in a strict initial segment 
   (or equivalently iff it is Lindel\"of), and
   unbounded otherwise. A function into $\LL_{\ge 0}$ or $\omega_1$ is bounded iff its range is bounded, and
   unbounded otherwise.
\end{defi}
When convenient, we see $\omega_1$ as a subset of $\LL_{\ge 0}$, identifying $\alpha\in\omega_1$
with $\langle\alpha,0\rangle\in\LL_{\ge 0}$. 
Hence, ``$\alpha\in\LL_{\ge 0}$'' should be understood as ``$\langle \alpha,0\rangle\in\LL_{\ge 0}$''.
The notions of stationary and club subsets of $\omega_1$ extend in an obvious way to subsets of $\LL_{\ge 0}$.
Recall the following classical lemma, whose proof can be found in any textbook about set theory.
\begin{lemma}[{\bf Fodor's Lemma}]
   Let $S\subset\omega_1$ be stationary and $f\co S\to\omega_1$ be (non necessarily continuous and)
   regressive, that is $f(\alpha)<\alpha$ $\forall\alpha\in S$. 
   Then there is $\alpha\in\omega_1$ such that $f^{-1}(\{\alpha\})$ is stationary.
\end{lemma}

\begin{defi}\label{def:general} 
   Let $Y$ be a space and $X\subset Y$ be a closed subspace.\\
   (a)
   $X$ is Type I in $Y$ iff there is a systematic cover $\{U_\alpha\,:\,\alpha\in\omega_1\}$ of $Y$
   such that $\wb{U_\alpha}\cap X$ is Lindel\"of for each $\alpha$. 
   When $X=Y$, we say that $X$ is Type I in itself, or just Type I.
   \\
   (b) A cover of $Y$ witnessing that $X$ is Type I in $Y$ is 
       a canonical cover of $X$ in $Y$ if it is good and systematic.\\
   (c)
   $X$ is functionally Type I in $Y$ (in short: f-Type I in $Y$) iff there is $s\co Y\to\LL_{\ge 0}$
   such that $X\cap s^{-1}([0,\alpha])$ is Lindel\"of  
   for each $\alpha\in\omega_1$.
   When $X=Y$, we say that $X$ is f-Type I in itself, or just f-Type I.
   In that case $s$ is called a slicer of $X$.
\end{defi}

The term {\em canonical} is motivated by the fact that two canonical covers
agree (on $X$) on a club set of indices, see Corollary \ref{cor:canonical} below.
The choice of $\LL_{\ge 0}$ for the range of $f$ in (c) (instead of another Type I space) is because
it plays a role similar to that of the interval $[0,1]$ in the classical Urysohn Lemma, or in
the definition of Tychonoff spaces. Indeed, if a space $X$ is Type I in itself and regular,
it is Tychonoff and f-Type I, with a slicer ``following'' a given canonical cover of $X$.
See Lemma \ref{lemma:TypeIregular} below for details.
Let us first clear up some trivialities.
\begin{lemma}\label{lemma:trivialTypeI} 
   Let $X$ be closed in a space $Y$.\\
   (a) 
   If $\{U_\alpha\,:\,\alpha\in\omega_1\}$ is a systematic cover of $Y$ 
   such that $\wb{U_\alpha}\cap X$ is Lindel\"of
   for each $\alpha<\beta$, then there is a canonical cover of $X$ in $Y$.\\
   (b) If $X$ is f-Type I in $Y$ then it is Type I in $Y$.\\
   (c) If $X$ is Type I in itself, then any closed subset of $X$ is 
   Type I in $X$. \\
   (d) If $X$ is f-Type I in itself, then any closed subset of $X$ is 
   f-Type I in $X$. \\
   (e) The Lindel\"of number of a space that is Type I in another is at most $\aleph_1$.
\end{lemma}
\proof \ \\
   (a)
   The only missing property is that $\cup_{\gamma<\alpha}U_\gamma = U_\alpha$ for limit $\alpha$.
   But in any chain cover we have $\wb{\cup_{\gamma<\alpha}U_\gamma}\subset \wb{U_{\alpha}}$, the latter being Lindel\"of,
   we may thus replace $U_\alpha$ by $\cup_{\gamma<\alpha}U_\gamma$, yielding a canonical cover of $X$.\\
   (b) If $s\co Y\to\LL_{\ge 0}$ witnesses that $X$ is f-Type I in $Y$, 
   then $\{s^{-1}([0,\alpha))\,:\,\alpha\in\omega_1\}$ is
   canonical for $X$ in $Y$. Indeed, $\wb{s^{-1}([0,\alpha))}$ is contained in $s^{-1}([0,\alpha])$ whose
   intersection with $X$ is Lindel\"of.\\
   (c), (d) \& (e) Immediate.
\endproof

\begin{defi} Let $Y$ be a space and $D\subset Y$ be a closed subspace.\\
   (a)
   We say that $D$ is narrow in $Y$ iff 
   given any
   systematic cover $\{U_\alpha\,:\,\alpha\in\omega_1\}$ of $Y$
   then either $D\subset U_\alpha$ for some $\alpha$ or 
   $\wb{U_\alpha}\cap D$ is Lindel\"of for each $\alpha$.
   \\
   (b)
   We say that $D$ is functionally narrow in $Y$ (in short: f-narrow in $Y$) iff 
   given any $f\co Y\to\LL_{\ge 0}$,
   either $f\upharpoonright D$ is bounded or $f^{-1}([0,\alpha])\cap D$ is Lindel\"of for each $\alpha$.
\end{defi}

\begin{rem}
   Slightly more general definitions arise by considering non-necessarily closed subsets and 
   asking for Lindel\"ofness of $\wb{U_\alpha}\cap \wb{D}$ and $f^{-1}([0,\alpha])\cap \wb{D}$,
   but we did not see any gain in choosing these versions.
\end{rem}

As above, we often abbreviate ``(f-)narrow in itself'' by just ``(f-)narrow''.
Let us immediately get rid of more trivialities. Point (c) is the main reason to 
restrict (f-)narrow subspaces to closed subsets.

\begin{lemma} Let $D$ be closed in a space $Y$.
   \\
   (a)
   If $D$ is narrow in $Y$, then $D$ is f-narrow in $Y$.\\
   (b) If $D$ is narrow (resp. f-narrow) in itself, then $D$ is narrow (resp. f-narrow) in $Y$.\\
   (c) If $C\subset D$ is closed and $D$ is narrow (resp. f-narrow) in $Y$, then $C$ is narrow (resp. f-narrow) in $Y$.
\end{lemma}
\proof\ \\
   (a) Given $s\co Y\to\LL_{\ge 0}$ which is unbounded on $D$,
    $\{s^{-1}([0,\alpha))\,:\,\alpha\in\omega_1\}$ is a systematic chain cover of $Y$.
   Then $\wb{s^{-1}([0,\alpha+1))}\cap D$ is Lindel\"of and contains
   $s^{-1}([0,\alpha])\cap D$, thus so is the latter.\\
   (b) \& (c) Immediate.
\endproof
Notice that the direction of the implication is the opposite between 
(f-)Type I and (f-)narrow spaces, as we have, in short (with the obvious warning about not writing explicitely the assumptions):
\begin{align*}
  \text{ f-Type I }&\Longrightarrow\text{ Type I }\\
  \text{ narrow }&\Longrightarrow\text{ f-narrow }.
\end{align*}

Also, notice that our definitions allow for trivial examples (take $U_\alpha$ to be the whole space for each
$\alpha\in\omega_1$):

\begin{example}\label{ex:Lindelof}
  Every closed Lindel\"of subspace of some space $Y$ is f-Type I and narrow, both in itself and in $Y$.
\end{example}

While we did not want to rule them out in order to avoid cumbersome exceptions in many proofs,
these are examples which we are not very interested in.

%%%%%%%%%%%%%%%%%%%%%%%%%%%%%%%%%%%%%%%%%%%%%%%%%%%%%%%%%%%%%%%%%%%%%%%%%%%%%%%%%%%%%%%%%%%%%%%%%%%%%%%%%%%%%%%%%%%%%%%%%%%%%%%%%%%

\section{Type I and functionally Type I subspaces}\label{sec:gen}

Let us start by stating an almost obvious useful lemma.
\begin{lemma}
   Let $\{U_\alpha\,:\alpha\in\kappa\}$ be a systematic chain cover of $Y$.
   Then a subset is closed [resp. open] iff its intersection with $\wb{U_\alpha}$ [resp. $U_\alpha$] is closed [resp. open]
   for each $\alpha$, and $f\co Y\to\LL_{\ge 0}$ is continuous iff its restriction to $U_\alpha$ is continuous for each $\alpha$
   iff its restriction to $\wb{U_\alpha}$ is continuous for each $\alpha$.
\end{lemma}
\proof
   Let $A\subset Y$.
   If $A$ is not closed, then there is $x\in\wb{A}-A$, choose $\alpha$
   such that $x\in U_\alpha$, then $A\cap \wb{U_\alpha}$ is not closed.
   If $A$ is not open then there is $x\in A-\text{int}(A)$, choose $\alpha$
   such that $x\in U_\alpha$, then $A\cap U_\alpha$ is not open.
   The rest follows immediately.
\endproof 

\begin{lemma}\label{lemma:Uknormal}
   Let $\mathcal{U} = \{U_\alpha\,:\,\alpha\in\omega_1\}$ be a systematic good cover of the space $Y$
   such that $U_\alpha\not= U_\beta$ whenever $\alpha\not=\beta$.
   If each $\wb{U_\alpha}$ is normal, then there is $s:Y\to\LL_{\ge 0}$ such that
   $s^{-1}([\alpha,\alpha+1])= \wb{U_{\alpha+1}} - U_\alpha$ 
   for each $\alpha\in\omega_1$. 
\end{lemma}
Notice that $s$ sends $\wb{U_{\alpha}} - U_\alpha$ (if nonempty) to $\{\alpha\}$ for each $\alpha\in\omega_1$.
\proof
   Since $A= \wb{U_{\alpha+1}} - U_{\alpha+1}$ and $B=\wb{U_\alpha}$ are disjoint 
   closed subsets of $\wb{U_{\alpha+1}}$, there is a Urysohn function 
   $s_\alpha\co \wb{U_{\alpha+1}}\to [\alpha,\alpha+1]\subset\LL_{\ge 0}$
   sending $A$ to $\alpha+1$ and $B$ to $\alpha$. (If $B$ is empty, just take a constant function
   on $\alpha+1$, if $A$ is empty but $B$ is not, a constant one on $\alpha$.)
   We may then define $s\co Y\to\LL_{\ge 0}$
   as equal to $s_\alpha$ on $\wb{U_{\alpha+1}} - U_{\alpha}$. Then $s$ is well defined, continuous,
   and satisfies the claimed properties. 
\endproof

As in Definition \ref{def:general} (c),
such an $s$ is called a {\em slicer of $\mathcal{U}$}, because it 
cuts the systematic cover into slices $s^{-1}(\{x\})$, respecting the inclusion relation
and sending the boundary of $U_\alpha$ to $\{\alpha\}$.
If $\mathcal{U}$ is a canonical cover of a Type I space $X$ and $s$ is a slicer of $\mathcal{U}$, 
we say that $s$ is a slicer of $X$ as well. 

\begin{lemma}\label{lemma:TypeIregular}
   Let $X$ be a closed subset of $Y$.\\ 
   (a) If $Y$ is regular and Type I (in itself) with canonical cover $\{V_\alpha\,:\,\alpha\in\omega_1\}$,
   then $Y$ is Tychonoff, each $\wb{V_\alpha}$ is normal, and  
   $Y$ is f-Type I (in itself).\\
   (b) If $Y$ is normal and $X$ is Type I in $Y$, then $X$ is f-Type I in $Y$.
\end{lemma}
\proof 
  If $X$ is Lindel\"of, there is nothing to do (see Example \ref{ex:Lindelof}). We thus assume that $X$ in non-Lindel\"of.
  \\
  (a)
  Almost immediate: a regular Lindel\"of space is normal, hence so is each $\wb{V_\alpha}$. 
  Then $s$ given by Lemma \ref{lemma:Uknormal} makes $Y$ f-Type I.
  Given $x$ in some open $U$, take $\alpha$ so that
  $x\in V_\alpha$, by normality of $\wb{V_\alpha}$ there is a real valued function that is $0$ on
  $x$ and $1$ on $\wb{V_\alpha} - (U\cap V_\alpha)$. Extend it to all of $Y$ by $1$ outside of $\wb{V_\alpha}$.
  This shows that $Y$ is Tychonoff.
  \\
  (b) By Lemma \ref{lemma:Uknormal}.
\endproof

It is trivial (and recorded in Lemma \ref{lemma:trivialTypeI} (c)--(d)) that 
(f-)Type I-ness (in a given space) is hereditary with respect to closed subspaces.
Let us now exhibit some examples that illustrate how this relative (f-)Type I-ness may
turn out in more general situations.
First, we note that
a non-Type I space might contain a closed subspace which is f-Type I in it.

\begin{example}
  $Y = \omega_2-\{\omega_1\}$ is not Type I (in itself) but $\omega_1\subset Y$
  is f-Type I in $Y$. 
  (More generally, the disjoint union of a f-Type I space and a space of Lindel\"of number $>\aleph_1$
   has the same property.)
\end{example}
\proof[Details] 
  It is obvious that $Y$ is not Type I (its Lindel\"of number is $\aleph_2>\aleph_1$).
  The function $f\co Y\to\omega_1$ which is the identity on $\omega_1$ and constant with any value on the rest
  shows that $\omega_1$ is f-Type I in $Y$.
\endproof

A space can become f-Type I with the removal of a single point.

\begin{example}
  The cone over $\LL_{\ge 0}$ defined as
  $Y = \LL_{\ge 0}\times[0,1] / \sim$, where $\langle x,1\rangle\sim \langle y,1\rangle$
  for each $x,y\in\omega_1$ is first countable and not Type I (in itself), but removing 
  its ``apex'' (the point with second coordinate $1$) yields an f-Type I space.
\end{example}
\proof[Details] 
  First countability follows from properties of the product $\LL_{\ge 0}\times[0,1]$:
  $\{\LL_{\ge 0}\times(1-1/n,1]\,:\,n\in\omega\}$ is a neighborhood base of $\LL_{\ge 0}\times\{1\}$,
  see Lemma \ref{lemma:propoct} (a) below.
  $Y$ is not Type I because no neighborhood of the apex point is Lindel\"of.
  The projection on the first factor 
  shows that $\LL_{\ge 0}\times[0,1)$ is f-Type I.
\endproof

A space can be Type I in itself but not in another space.

\begin{example}\label{ex:TP}
  The closed 
  subspace $\omega_1\times\{\omega\}$ is not Type I in the Tychonoff plank 
  $T = (\omega_1+1)\times(\omega+1) - \{\langle\omega_1,\omega\rangle\}$
  but is f-Type I in itself. 
\end{example}
\proof[Details]
   $(\omega_1+1)\times\omega$ is Lindel\"of (even $\sigma$-compact) and dense in $T$, hence any 
   systematic chain cover of $T$ must stagnate (and contain the whole space) after at most countably many steps.
\endproof

One may believe that 
having a normal ambient space does prevent this situation, but it is not the case.
\begin{example}\label{ex:gammaN}
  The copy of $\omega_1$ in any $\gamma\N$ 
  is not Type I in $\gamma\N$, 
  but is obviously f-Type I in itself. Moreover, $\gamma\N$ is a normal space.
\end{example}
\proof[Details]
  $\gamma\N$ is a symbol used in particular by P. Nyikos to denote
  a normal space that is the union of a countable dense subset of isolated points
  (identified with the integers) and a closed copy of $\omega_1$, see e.g. \cite{Nyikos:2003}.
  Such spaces exist in {\bf ZFC}.
  Any systematic chain cover of $\gamma\N$ stagnates after at most countably many steps.
\endproof

We now turn to the problem of when a (relative) Type I space is (relatively) f-Type I.
The simplest instance is for the ambient space itself: when does Type I in itself
imply f-Type I in itself.
Lemma \ref{lemma:TypeIregular} shows that this is the case for regular spaces. 
We now present two examples of Hausdorff spaces for which the implication fails.
The first one is given by ``piling up'' copies of Roy's lattice space (see e.g. \cite[Ex. 126]{CEIT}),
whose real-valued maps are constant, preventing f-Type I-ness.
The second one is even {\em functionally Hausdorff} (also called {\em completely Hausdorff}), that is,
given $x,y\in X$, there is a 
real valued fonction on $X$ with $f(x) = 0$, $f(y) = 1$, but its construction 
is a bit more cumbersome. 

\begin{example}\label{ex:conHaus}
  A first countable connected Hausdorff Type I space $R$ which is not
  f-Type I.
\end{example}
(Note: This is a corrected version of the space constructed in \cite[Lemma 10.6]{mesziguesDirections} whose description is flawed. )
\proof[Details]
As a set, $R$ is a subset of $\LL_{\ge 0}\times{(\omega+1)}$.
Let $\{C_n\,:\,n\in\omega\}$ be a disjoint collection of dense subsets of $\Q\cap[0,1)$. We assume that $0\in C_0$.
We then let $L_n$ to be the subspace of $\LL_{\ge 0}$ defined as $\omega_1\times C_n$ (with lexicographic order topology).
Define $R$ as $\cup_{n\in\omega} L_n\times\{n\} \cup \omega_1\times\{\omega\}$,
with the following topology.
Denote by $q_\alpha$ the point $\langle\alpha,\omega\rangle\in R$. Given $a<b\in\LL_{\ge 0}$ and $n<\omega$, denote by 
$(a,b)_n$
the set $L_n\cap (a,b)\subset \LL_{\ge 0}$ and 
set
$U_{a,b,n} = (a,b)_n\times\{n\}$. Let $z = \langle x,n\rangle\in R$, with $n<\omega$.
A neighborhood of $z$ is given by $U_{a,b,n}$ 
if $n$ is even and by a stack of $3$ intervals $U_{a,b,n-1}\cup U_{a,b,n}\cup U_{a,b,n+1}$
if $n$ is odd, with $a<x<b$ in each case.
Finally, a neighborhood of $q_\alpha$
is $\{q_\alpha\}\cup\cup_{n\ge m} (\alpha,\alpha+1)_n\times\{n\}$, for $m$ even.
Set $R_\alpha = R \cap [0,\alpha)\times(\omega+1)$.

\begin{claim}\label{claim:RTypeI}
  $R$ is Hausdorff, Type I and each $R_\alpha$ is connected (and thus $R$ too).
\end{claim}
\proof
  Hausdorffness is immediate since the $C_n$ are disjoint.
  By construction, since $0\in C_0$,
  $\wb{R_\alpha}=R_{\alpha}\cup\{\langle\alpha,0\rangle\}$ ($\alpha$ is seen as a member of $\LL_{\ge 0}$ here) 
  is countable and thus Lindel\"of, so $R$ is Type I.
  Now, $R_{\alpha+1}-R_\alpha$ 
  is exactly Roy's lattice space and is thus connected (see details in \cite[Ex. 126]{CEIT}).
  Since $R_{\alpha+1} = (R_{\alpha+1}-R_\alpha)\cup\wb{R_\alpha}$,
  $(R_{\alpha+1}-R_\alpha)\cap\wb{R_\alpha}= \{\langle\alpha,0\rangle\}\not=\varnothing$,
  and $R_\alpha = \cup_{\beta<\alpha}R_\alpha$ when $\alpha$ is limit,
  it follows by induction that $R_\alpha$ is connected (and thus so is $\wb{R_\alpha}$).
\endproof

\begin{claim}\label{claim:Rnottrue}
  Any $g:R \to \LL_{\ge 0}$ is constant.
\end{claim}
\proof
  Since $R_\alpha$ is connected and countable,
  $g$
  must send $R_\alpha$ to a connected and at most countable subspace of $\LL_{\ge 0}$.
  The only such subspaces of $\LL_{\ge 0}$ are singletons, hence
  $g$ must be constant on $R_\alpha$ and thus on all of $R$.  
\endproof

\noindent
It follows at once from Claim \ref{claim:Rnottrue} that $R$ cannot be f-Type I.
\endproof

Our next example is a subset of a (set-theoretic) {\em tree}.
Recall that a 
tree is a partially ordered set such that the set of predecessors of each member is well ordered.
Our terminology is standard, see e.g. \cite{Nyikos:trees} if in need of definitions of the 
terms {\em height, level, chain} and {\em antichains}.

\begin{figure}
  \begin{center}
    \epsfig{figure = 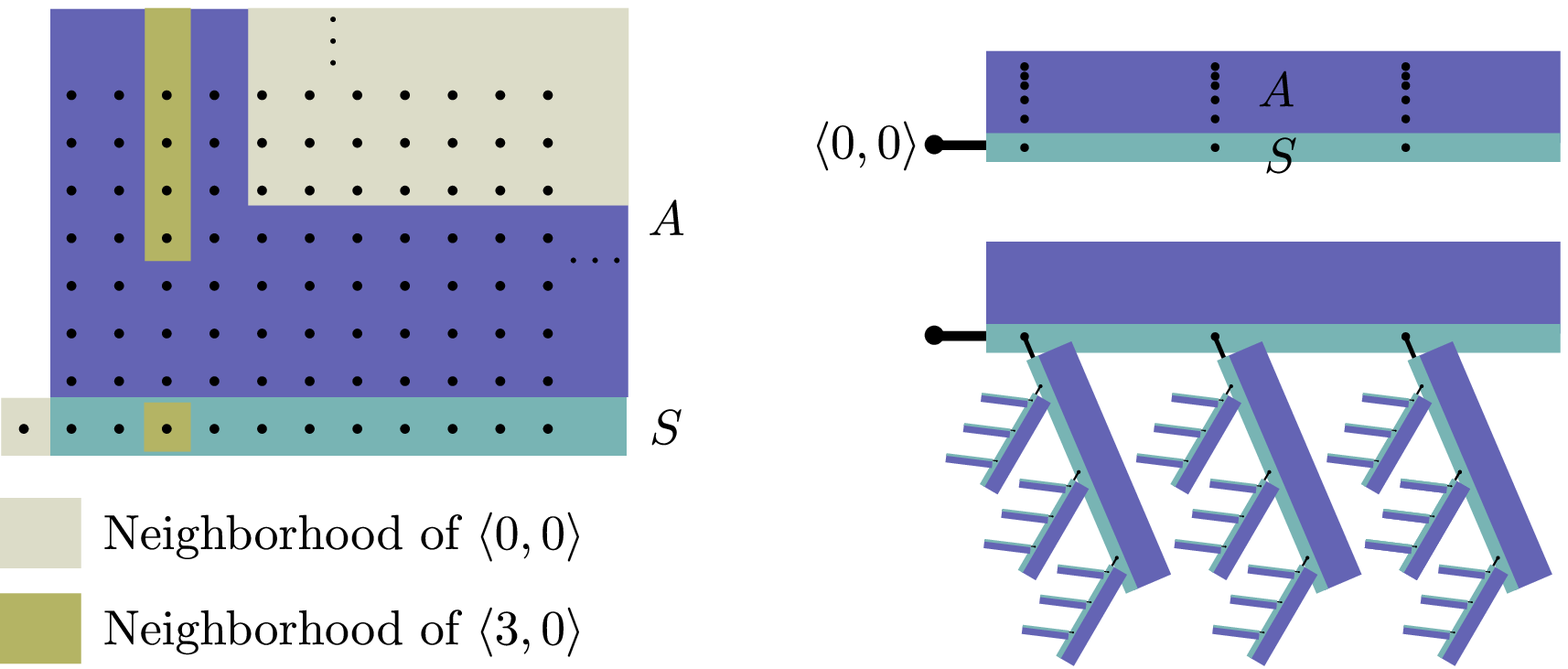, width = .85\textwidth}
    \caption{Example \ref{ex:confuncHaus}}
    \label{fig:Roy}
  \end{center}
\end{figure}

\begin{example}\label{ex:confuncHaus}
  A first countable functionally Hausdorff Type I space $T$ which is not
  f-Type I.
\end{example}
\proof[Details]
  $T$ is a subset of a tree homeomorphic to ${}^{<\omega_1}\omega$, 
  but its topology is weaker than the order topology. 
  We first recall an example of a functionally Hausdorff non-regular countable space \cite[Ex. 79]{CEIT}
  which we will use as a building brick (see Figure \ref{fig:Roy}, left).
  (Beware that functionally Hausdorff spaces are called {\em Urysohn},
  and completely Hausdorff means another closely related property, in \cite{CEIT}.\footnote{Separation axioms are a bit like
  the {\em moiti\'e-moiti\'e}, which may be understood as a dish of melted cheese or an alcoholic beverage
  depending on whether you are in Geneva or in Fribourg, actually even depending on the establishment on which
  you are sitting in Geneva. But we digress.})
  Let $B$ be the subset of the lattice of integers $\omega\times\omega$ defined as
  $\{\langle 0,0\rangle\}\cup A\cup S$, where $A=\{\langle n,m\rangle\,:\,n,m\in\omega,\, n,m\ge 1\}$ 
  and $S = \{\langle n,0\rangle\,:\,n\in\omega,\, n\ge 1\}$.
  Points in $A$ are isolated, while a neighborhood base for $\langle n,0\rangle\in S$
  is given by $\{\langle n,0\rangle\}\cup\{\langle n,m\rangle\,:\,m\in\omega, m > k\}$ for $k\in\omega$.
  That is, points in the vertical line above $\langle n,0\rangle$ converge to it.
  A neighborhood base of $\langle 0,0\rangle$ is given by 
  $\{\langle 0,0\rangle\}\cup\{\langle n,m\rangle\,:\,n,m\in\omega, n,m > k\}$ for $k\in\omega$.
  It is easy to see that $S$ is functionally Hausdorff but non-regular 
  since $S$ is closed in it and cannot be separated from $\langle 0,0\rangle$ (see \cite[Ex. 79]{CEIT} for details).
  In fact, any neighborhood of $\{\langle 0,0\rangle\}$ has all the points in $S$ with first coordinate $>k$ (for some $k$)
  in its closure.
  Set $B_0$ to be $B-\{\langle 0,0\rangle\} = A\cup S$.
  \\
  We now define $T$, which is a subset of ${}^{<\omega_1}B_0$.
  We see members
  of $T$ as sequences of length $<\omega_1$ with entries in $B_0$.
  Recall that the order $\le$ on ${}^{<\omega_1}B_0$ is given by sequence extension,
  and that $\sigma^\smallfrown b$ denotes the extention of the sequence $\sigma\in {}^{<\omega_1}B_0$
  by adding $b\in B_0$ as the last member. 
  We define $T$ (and its topology) by induction on the levels. 
  As one may guess,
  $T_\alpha$, $T_{<\alpha}$ and $T_{\le\alpha}$ respectively mean the elements at level $\alpha$,
  $<\alpha$ and $\le\alpha$. A canonical cover for $T$ will be given
  by $\{T_{<\alpha}\,:\,\alpha\in\omega_1,\,\alpha\text{ is limit}\}$.
  Each $\sigma\in T_{\le\alpha}$ will have a countable local base in $T_{\le\alpha}$
  denoted by $\mathcal{B}(\sigma,\alpha)$.\\
  The idea of the construction is the following. We start with a copy of $B$, and at successor levels
  glue copies of $B$ to each member $s$ of $S$, identifying $s$ with $\langle 0,0 \rangle$
  in the glued copy. The first four stages of this construction are
  illustrated in Figure \ref{fig:Roy}, right.
  At limit levels, since we are in a subtree of a tree homeomorphic to ${}^{<\omega_1}\omega$,
  we take the topology of cylinders (actually, equivalent to it)
  to obtain (at this level) a space homeomorphic to $\R-\Q$. 
  $S$ being a closed discrete subspace of $B$ and the copies of $B_0$ being discrete as well
  in $T$ at successor levels, we cannot have uncountably many of them 
  if we want the space to be of Type I. Hence, after a limit level, we extend only countably
  many sequences, choosing a countable subset which is dense in this limit level.
  \\
  We now describe the construction in more details. The set $C_\alpha\subset T_\alpha$ will contain
  the sequences that continue in the next levels.
  $T_0 = C_0$ contains the empty sequence, which is open in $T_0$.
  The $\alpha+1$th level $T_{\alpha+1}$ is given by the set of $\sigma^\smallfrown b$ for each $b\in B_0$ and each
  $\sigma\in C_{\alpha}$, and $C_{\alpha+1}$ is the subset
  of those $\sigma^\smallfrown b$ with $b\in S$. Then $\{\sigma\}\cup\{\sigma^\smallfrown b\,:\,b\in B_0\}$
  will be a copy of $B$, with $\sigma$ playing the role of $\langle 0,0\rangle$.
  The local bases (in $T_{\le\alpha+1}$) of points at levels $<\alpha$ do not change.
  Given $\sigma\in C_\alpha$ and $O\in\mathcal{B}(\sigma,\alpha)$, set
  $$ V_k(O) = O\cup\{\mu^\smallfrown \langle n,m\rangle\,:\, \mu\in O\cap C_\alpha,\,n,m > k\}.$$
  Then $\mathcal{B}(\sigma,\alpha+1) = \{V_k(O)\,:\,O\in\mathcal{B}(\sigma,\alpha),\,k\in\omega\}$.
  For $\langle n,0\rangle\in S$, 
  $\mathcal{B}(\sigma^\smallfrown \langle n,0\rangle,\alpha+1)$ is the set containing
  $\{\sigma^\smallfrown \langle n,0\rangle\}\cup\{\sigma^\smallfrown\langle n,m\rangle\,:\, m > k\}$
  for each $k\in\omega$.
  Finally, if $a\in A$, then $\sigma^\smallfrown a$ is isolated.
  \\
  If $\alpha$ is limit, we set $T_\alpha$ to be 
  $\{ \sigma\in{}^\alpha B_0\,:\,\sigma\upharpoonright\beta\in C_\beta\,\forall\beta<\alpha\}$.
  The neighborhood bases in $T_{\le\alpha}$ do not change for points strictly below level $\alpha$.
  We add to the topology the bases for points in
  $T_\alpha\subset T_{\ge\alpha}$ given by 
  $\{ U_\sigma\,:\,\sigma\in \cup_{\beta < \alpha} C_\beta\}$, where 
  $U_\sigma = \{\mu\in T_{\le\alpha}\,:\,\mu>\sigma\}$.
  Notice that $U_\sigma$ and $U_\nu$ are either disjoint (when $\sigma$ and $\nu$
  are incomparable) or included one in the other (when $\sigma$ and $\nu$
  are comparable), and that $U_\sigma$ is open in $T_{<\alpha}$.
  This topology is stronger than the usual cylinder topology on $T_\alpha$, where one fixes only finitely many values,
  but still secound countable. Actually, we have:
  \begin{claim}
      $T_\alpha$ is homeomorphic to $\R-\Q$ when $\alpha$ is limit.
  \end{claim}
  \proof
      Let $\alpha_n$ be an increasing sequence of ordinals with supremum $\alpha$.
      Then $T_\alpha$ is homeomorphic to the space $M$ of $\omega$-sequences whose $n$th 
      members are is $C_{\alpha_n}$ with the initial segment topology -- that is, open
      sets are defined by fixing the first entries of the sequence. (Just forget the entries of 
      the $\alpha$-sequence $\sigma\in T_\alpha$ that are in between levels $\alpha_n$.)
      By construction, each $C_{\alpha_n}$ is countable and given $\mu\in C_{\alpha_n}$,
      there are $\omega$-many members of $C_{\alpha_{n+1}}$ above $\mu$.
      Hence, $M$ is homeomorphic to ${}^\omega\omega$ with the initial segment topology, which
      in turn is homeomorphic to $\R-\Q$.
  \endproof
  We then let $C_\alpha$ be a countable dense subset of $T_\alpha$.
  This defines $T$ and its topology. By construction, each $T_{\le\alpha}$ is second countable.
  \begin{claim}
      $\mathcal{U} = \{T_{<\alpha}\,:\,\alpha\text{ is limit}\}$ is a canonical cover
      making $T$ a Type I space.
  \end{claim}
  \proof
     By construction, $\wb{T_{<\alpha}} = T_{\le\alpha}$ when $\alpha$ is limit.
  \endproof
  Notice that $T_{<\alpha+1}=T_{\le\alpha}$ is not open in $T$, but 
  $T_{\le\alpha}\cup\{\mu^\smallfrown\langle n,m\rangle\,:\,\mu\in C_{\alpha},\,n,m>0\}$ is.
  \begin{claim}
      $T$ is functionally Hausdorff.
  \end{claim}
  \proof
      Notice that any real valued function $f$ on $T_{\le\alpha}$ can be continuously extended to all of $T$
      by letting $f(\sigma) = f(\sigma\upharpoonright\alpha)$ for $\sigma$ at a level higer than $\alpha$.
      Let $\sigma,\mu$ be distinct members of $T$. If one of them is isolated, we are over.
      Otherwise, none of their entries (as sequences) is in $A$.
      If $\sigma$ and $\mu$ are incomparable in $^{<\omega_1}B_0$, let $\alpha$ be minimal such that
      $\sigma\upharpoonright\alpha\not=\mu\upharpoonright\alpha$. 
      Then $\alpha=\beta+1$ for some $\beta$, 
      hence $\sigma\upharpoonright\alpha = \nu^\smallfrown \langle n,0\rangle$ and 
      $\mu\upharpoonright\alpha = \nu^\smallfrown \langle m,0\rangle$ for some $\nu$ at level $\beta$,
      and $n\not=m$. Define $f\co T_{\le\alpha}\to\R$ to take value
      $1$ on $\{\nu^\smallfrown \langle n,k\rangle\,:\,k\in\omega\}$
      and $0$ elsewhere. Then $f$ is continuous, extending it above $T_{\le\alpha}$
      yields a map $f\co T\to \R$ with $f(\sigma) = 1$, $f(\mu)=0$.
      If $\sigma<\mu$, let $n$ be such that $\sigma^\smallfrown \langle n,0\rangle\le\mu$.
      Let $f$ take value $1$ on  
      $\{\sigma^\smallfrown \langle n,k\rangle\,:\,k\in\omega\}$
      and all the points above $\sigma^\smallfrown \langle n,0\rangle$, and $0$ elsewhere.
      Then $f$ is continuous, and takes value $0$ on $\sigma$ and $1$ on $\mu$.
  \endproof
  \begin{claim}
      If $\sigma\in C_\alpha\subset T$ and $f\co T\to\R$, then for any interval $(a,b)\ni f(\sigma)$
      and any $\beta \ge \alpha$, there is $\mu\in C_\beta$, $\mu \ge \sigma$, such that $f(\mu)\in[a,b]$.
  \end{claim}
  \proof
      By induction on $\beta$.
      If $\beta=\alpha$, it is obvious. 
      If $\beta = \gamma+1$, fix $c,d\in\R$ such that $f(\sigma)\in(c,d)\subset[c,d]\subset(a,b)$ and 
      $\mu\in f^{-1}([c,d])\cap C_\gamma$, $\mu\ge \sigma$.
      By definition of the neighborhood base, $\wb{f^{-1}((a,b))}$ contains 
      $\mu^\smallfrown \langle n,0\rangle$ if $n$ is big enough, hence so does $f^{-1}([a,b])$.
      Assume now that $\beta$ is limit and let $\beta_n$ be an increasing sequence converging to $\beta$.
      Let $d$ be the minimum between $f(\sigma)-a$, $b-f(\sigma)$ and let $d_n$ be a sequence of 
      strictly positive numbers whose sum is less than $d/2$.
      By induction, for each $n$ there is $\nu_n\in C_{\beta_n}$
      such that $\nu_{n+1} \ge \nu_n$ (hence $\nu_{n+1}$ extends $\nu_n$)
      and $f(\nu_{n+1})\in [f(\nu_n)-d_n,f(\nu_n)+d_n]$.
      By construction, $\nu$ defined as the supremum (in $T$) of the $\nu_n$ (which thus converge to $\nu$)
      is a member of $T_\beta$, and $|f(\sigma)-f(\nu)|\le d/2$.
      By density (and the fact that the members of $T_\beta$ above $\sigma$ form an open set),
      there is a member $\mu$ of $C_\beta$ above $\sigma$ in $f^{-1}((f(\nu)-d/2,f(\nu)+d/2))$.
      Since $|f(\mu)-f(\sigma)|\le d$, $f(\mu)\in[a,b]$.
  \endproof
  If follows that $T$ is not f-Type I. Indeed, given $f\co T\to \LL_{\ge 0}$,
  take any $\sigma\in C_\alpha\subset T$ and $\gamma$ such that $f(\sigma)\in[0,\gamma]$.
  Define $g\co T\to [0,\gamma+2]$ as the minimum between $\gamma+2$ and $f$,
  then $f^{-1}([0,\gamma+1]) = g^{-1}([0,\gamma+1])$ is unbounded in the tree 
  by the previous claim and hence non-Lindel\"of.
\endproof

\begin{rem} 
   The space $B$ has a weaker regular topology, but 
   there is no weaker regular topology on $T$ such that $\{T_{<\alpha}\,:\,\alpha\in\omega_1\}$ is still a canonical cover.
   Indeed,
   if $X$ with topology $\tau$ is Type I with canonical cover 
   $\{U_\alpha\,:\,\alpha\in\omega_1\}$ 
   and there is a weaker topology $\rho\subset\tau$ such that $\langle X,\rho\rangle$
   is regular, $U_\alpha$ is $\rho$-open for each $\alpha$,
   and the $\rho$-closure of $U_\alpha$ is included in $U_\gamma$ for some $\gamma>\alpha$,
   then by Lemma \ref{lemma:TypeIregular}, $\langle X,\tau\rangle$ is f-Type I. 
\end{rem}

By Lemma \ref{lemma:TypeIregular} (b),  
being relatively Type I is equivalent to being relatively f-Type I when the ambient space is normal.
We were not able to decide the following question, although we anticipate a positive answer.
 \begin{q} 
    Is there a regular (or even Tychonoff) non-Type I space $Y$ containing a subspace $X$ which is Type I in $Y$ 
    but not f-Type I in $Y$~?
\end{q}

\vskip .3cm

Let us now briefly mention a partial order on systematic covers of spaces,
whose study we have not pursued (except for the simple lemmas below), 
but which gives another light on the concept of a canonical cover.
If $\mathcal{U}$ is a cover of $X$ and $B\subset X$, we let 
$\mathcal{U}\cap B$ be $\{U\cap B\,:\,U\in \mathcal{U}\}$.
\begin{defi}
   Let $\mathcal{U},\mathcal{V}$ be covers of some space $X$, and $B\subset X$ be closed.
   Then $\mathcal{U}\trianglelefteq_B\mathcal{V}$ iff 
   $\mathcal{U}\cap B$ is a refinement of 
   $\mathcal{V}\cap B$, that is,
   for each $U\in\mathcal{U}$ there is $V\in\mathcal{V}$
   with $U\cap B\subset V\cap B$. 
   If $\mathcal{U}\trianglelefteq_B\mathcal{V}\trianglelefteq_B\mathcal{U}$ we write $\mathcal{U}\bowtie_B\mathcal{V}$.
   We abbreviate $\trianglelefteq_X,\bowtie_X$ by $\trianglelefteq,\bowtie$.
\end{defi}
Our first lemma shows that canonical covers are in a sense $\trianglelefteq$-minimal.
\begin{lemma}
   \label{lemma:canomin}
   Let $B\subset X$ be Type I in $X$,
   $\mathcal{U}$ be a canonical cover of $B$ in $X$
   and $\mathcal{V}=\{V_\alpha\,:\,\alpha\in\omega_1\}$ be a chain cover of $X$.
   Then $\mathcal{U}\trianglelefteq_B \mathcal{V}$.
\end{lemma}
\proof
   The closure of each member of $\mathcal{U}$ intersected with $B$ is Lindel\"of, and hence
   included in some $V_\alpha$.
\endproof

\begin{lemma}\label{lemma:U_alphaV_alpha}
   Let $\mathcal{U}=\{U_\alpha\,:\,\alpha\in\lambda\},\mathcal{V}=\{V_\alpha\,:\,\alpha\in\kappa\}$ 
   be good systematic covers of some space and $B$ be a closed subset. Assume
   that $U_\alpha\not= U_\beta$ and $V_\alpha\not= V_\beta$ when $\alpha<\beta$.
   Then $\mathcal{U}\bowtie_B\mathcal{V}$ iff $\lambda = \kappa$ and 
   $\{\alpha\in\kappa\,:\, U_\alpha\cap B = V_\alpha\cap B\}$ is club in $\kappa$.
\end{lemma}
\proof
The reverse implication is immediate.
For the direct implication,
recall that $\lambda,\kappa$ are regular cardinals since the covers are good.
For each $\alpha\in\lambda$ there are $\beta\in\kappa$ and $\gamma\in\kappa$ such that 
$$U_\alpha\cap B\subset V_{\beta}\cap B\subset \wb{V_{\beta}}\cap B\subset U_{\gamma}\cap B,$$ 
hence
$\kappa$ and $\lambda$ have the same cofinality and are thus equal.
A leapfrog argument then shows that $\{\alpha\in\kappa\,:\, U_\alpha\cap B = V_\alpha\cap B\}$ is club.
\endproof

The following well known corollary explains the choice of terminology for canonical covers.
(This terminology comes from P. Nyikos, who uses the slightly different term {\em canonical sequence} in \cite{Nyikos:1984}.)

\begin{cor}\label{cor:canonical} 
   Let $B$ be closed and Type I in the space $X$.
   If $\mathcal{U}=\{U_\alpha\,:\,\alpha\in\omega_1\},\mathcal{V}=\{V_\alpha\,:\,\alpha\in\omega_1\}$ 
   are canonical covers of $B$ in $X$, then 
   $\{\alpha\,:\, U_\alpha\cap B = V_\alpha\cap B\}$ is club in $\omega_1$.
\end{cor}
\proof
   Lemma \ref{lemma:canomin} implies $\mathcal{U}\bowtie_B\mathcal{V}$, and we conclude with Lemma \ref{lemma:U_alphaV_alpha}.
\endproof

We finish this section with some definitions and results that will be useful later.
If $\mathcal{U}=\{U_\alpha\,:\,\alpha\in\omega_1\}$ is a canonical cover of a Type I subspace $D$ of $X$,
its set of {\em bones} is $\BO(\mathcal{U},D)=\{(\wb{U_\alpha}-U_\alpha)\cap D\,:\,\alpha\in\omega_1\}$, and its {\em skeleton} 
$\text{Sk}(\mathcal{U},D)$ is the union of the bones. Notice that $\text{Sk}(\mathcal{U},D)$ is closed in $X$.
We denote $\BO(\mathcal{U},X)$, $\text{Sk}(\mathcal{U},X)$ by $\BO(\mathcal{U})$,
$\text{Sk}(\mathcal{U})$. (This terminology is also due to P. Nyikos.)
Recall that a space is {\em $\omega_1$-compact} (or has {\em countable extent}) iff
its closed discrete subspaces are at most countable.
\begin{lemma}\label{lemma:skeleton}
   Let $D\subset X\subset Y$, where $D$ is closed non-Lindel\"of, $X$ closed and Type I in the space $Y$.
   Let $\mathcal{U}$ be a canonical cover of $X$ in $Y$.
   Then the following hold.\\
   (a) If $X$ is countably compact, then $D$ intersects the members of $\BO(\mathcal{U})$ on a club set of indices.\\
   (b) If $X$ is $\omega_1$-compact, then $D$ intersects the members of $\BO(\mathcal{U})$ on a stationary set of indices.
\end{lemma}
\proof 
   This is well known, and (a) is straightforward. For (b), assume that $X$ is $\omega_1$-compact.
   If $D$ misses the bones of $\mathcal{U}$ on a club set of indices,
       taking one point of $D$ (when available) between each member of this club set 
       yields a closed discrete uncountable subset of $X$.
\endproof

Recall that a space is {\em $\omega$-bounded} iff each countable subset has a compact closure.
A product of $\omega$-bounded spaces is $\omega$-bounded.
The following lemma is immediate.
\begin{lemma}\label{lemma:omegabounded}
   Let $X$ be closed and Type I in the space $Y$. Then $X$ is countably compact iff $X$ is $\omega$-bounded
   iff each member of a canonical cover of $X$ has compact closure when intersected with $X$.
\end{lemma}

%%%%%%%%%%%%%%%%%%%%%%%%%%%%%%%%%%%%%%%%%%%%%%%%%%%%%%%%%%%%%%%%%%%%%%%%%%%%%%%%%%%%%%%%%%%%%%%%%%%%%%%%%%%%%%%%%%%%%%%%%%%%%%%%%%%

\section{Narrow and functionally narrow subspaces}\label{sec:narrow}

Notice that the Tychonoff plank (Example \ref{ex:TP}) and 
any $\gamma\N$ (Example \ref{ex:gammaN}) are narrow in themselves but not Type I in themselves.
We will be more interested in cases where narrow (sub)spaces are also Type I.
We start with the most basic examples.
Recall that the longline $\LL$ consists of two copies of $\LL_{\ge 0}$ glued at their $0$ point, 
with the reverse order on one of the copies. 

\begin{example}\label{ex:1st}
   $\LL_{\ge 0}$, $\omega_1$ are Type I and narrow in themselves.
   $\LL$, $(\LL_{\ge 0})^2$ are Type I in themselves but not f-narrow.
\end{example}
\proof[Details] Anyone who read until this point should be convinced, or see the next result.
\endproof

The example of $\LL$ shows in passing that the union of two closed subsets which are (f-)narrow in some space
is not always (f-)narrow.
It might be interesting to note
the following (easy) theorem.

\begin{thm}\label{thm:stationary}
   Let $S$ be an uncountable subset of $\omega_1$ endowed with the subspace topology.
   Then $S$ is f-Type I and the following are equivalent.\\
   (a) $S$ is stationary,\\
   (b) $S$ is $\omega_1$-compact,\\
   (c) $S$ is narrow in itself.
\end{thm}
\proof
   It is immediate that $S$ is f-Type I, and
   (a) $\leftrightarrow$ (b) is classical. 
   \\
   (a) $\rightarrow$ (c). If $A$ is closed non-Lindel\"of in $S$, it is unbounded.
   Hence $A$ contains $C\cap S$ for a club $C\subset\omega_1$, hence a stationary subset of $\omega_1$.
   By Fodor's lemma, any open $V\supset A$ contains a terminal part of $S$.
   Since the rest of the space is countable, any systematic (actually, chain) cover of $S$ with a non-Lindel\"of
   member has a member containing all of $S$.\\
   (c) $\rightarrow$ (a). If $S$ is non-stationary, it misses a club $C\subset\omega_1$.
   Take one isolated point of $S$ between each member of $C$ (when available), this yields
   a clopen uncountable discrete subset $E$ of $S$. By throwing away some points, we may assume that $E$ has
   uncountable complement in $S$.
   Then $\{([0,\alpha)\cap S)\cup E\,:\,\alpha\in\omega_1 \}$ is a systematic cover of $S$
   whose members have non-Lindel\"of closure. 
\endproof

\begin{figure}
  \begin{center}
    \epsfig{figure = 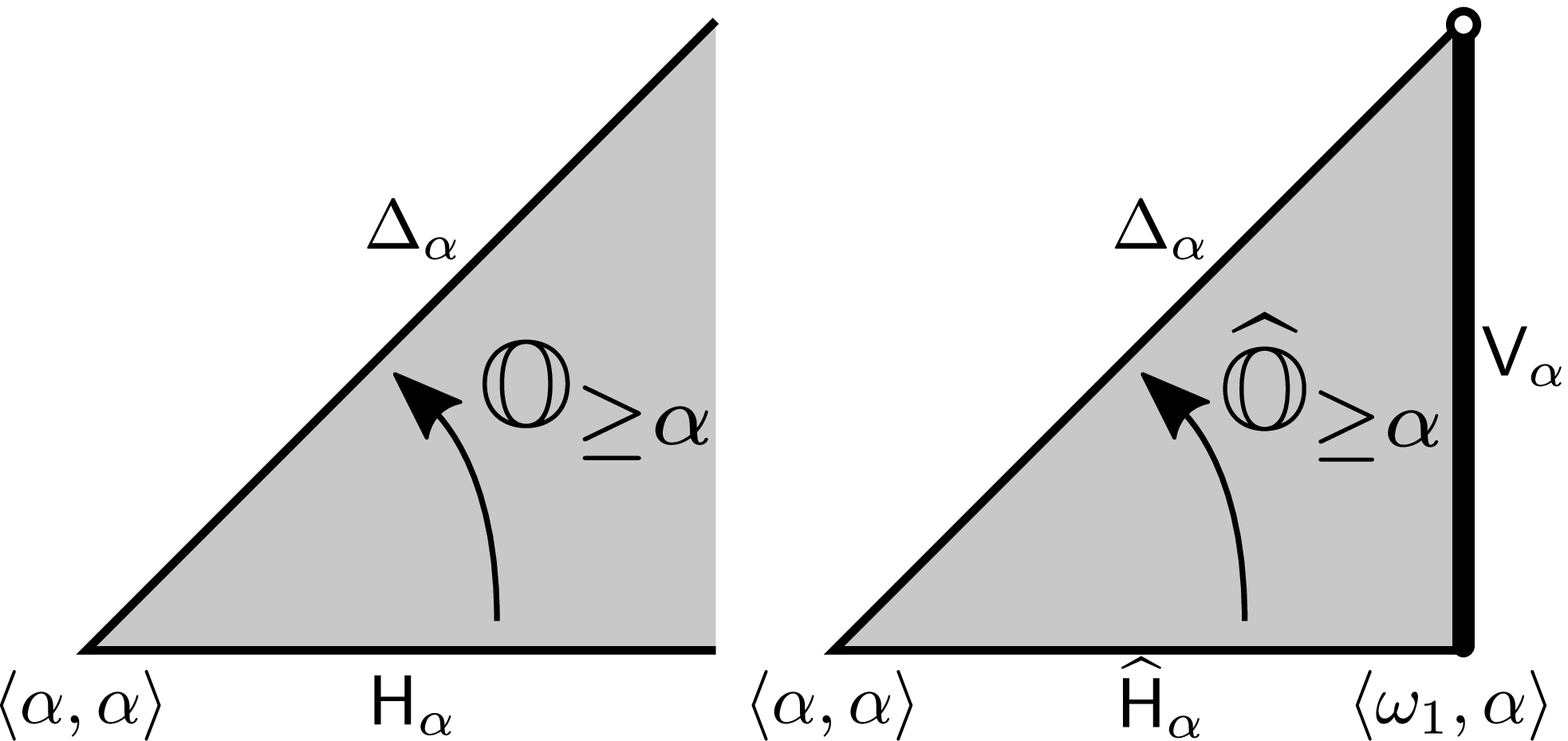, width = .5\textwidth}
    \caption{Octants}
    \label{fig:octants}
  \end{center}
\end{figure}

We now define two spaces which we use as building blocks in many of our examples.
When considering the space $\LL_{\ge 0}\cup\{\omega_1\}$, the neighborhoods of $\{\omega_1\}$
are of course final segments of $\LL_{\ge 0}$ union $\{\omega_1\}$.
\begin{defi}\ \\ 
   The
   octant is defined as $\mathbb{O} = \{\langle x,y\rangle\in(\LL_{\ge 0})^2\,:\,y\le x\}$,
   and its diagonal and horizontal ``boundaries'' are denoted $\Delta =\{\langle x,x\rangle\,:\,x\in\LL_{\ge 0}\}$,
   $\HH = \{\langle x,0\rangle\,:\,x\in\LL_{\ge 0}\}$.
   \\
   The ``horizontally compactified octant'' 
   is $\wh{\mathbb{O}} = \{\langle x,y\rangle\in\LL_{\ge 0}\cup\{\omega_1\}\times\LL_{\ge 0}\,:\,y\le x\}$,
   with its diagonal, horizontal and vertical ``boundaries''  
   $\Delta$, $\wh{\HH}=\{\langle x,0\rangle\,:\,x\in\LL_{\ge 0}\cup\{\omega_1\}\}$,
   $\VV = \{\langle \omega_1,x\rangle\,:\,x\in\LL_{\ge 0}\}$. \\
   We also define $\mathbb{O}_{\ge\alpha} = \{\langle x , y\rangle\in\mathbb{O}\,:\,x,y\ge\alpha\}$
   (and $\Delta_\alpha,\HH_\alpha$)
   and $\wh{\mathbb{O}}_{\ge\alpha} = \{\langle x , y\rangle\in\wh{\mathbb{O}}\,:\,x,y\ge\alpha\}$
   (and $\Delta_\alpha, \wh{\HH}_\alpha, \VV_\alpha$ defined accordingly).
\end{defi}
$\mathbb{O}_{\ge\alpha}$, $\wh{\mathbb{O}}_{\ge\alpha}$ are represented as in Figure \ref{fig:octants}.
$\mathbb{O}_{\ge\alpha}$ is homeomorphic to $\mathbb{O}_{\ge\beta}$ for each $\alpha,\beta\in\omega_1$, 
and so are $\wh{\mathbb{O}}_{\ge\alpha}$ and $\wh{\mathbb{O}}_{\ge\beta}$.
Notice that $\mathbb{O}_{\ge\alpha}$ is first countable while 
$\wh{\mathbb{O}}_{\ge\alpha}$ is not at the points of $\VV_\alpha$,
hence the bolder line used for representing $V_\alpha$ in Figure \ref{fig:octants}.
The arrow from $H_\alpha$ to $\Delta_\alpha$ is inspired by Point (c) in the next well known lemma, which
contains most of the basic topological properties of these octants. 
The notions ``stationary'' and ``club'' extend naturally to copies of $\omega_1$ or $\LL_{\ge 0}$
in $(\LL_{\ge 0})^2$ such as $\Delta,\mathsf{H},\mathsf{V}$.

\begin{lemma}[{\bf Properties of $(\LL_{\ge 0})^2,\mathbb{O},\wh{\mathbb{O}}$}]
   \label{lemma:propoct}
   \ \\
   (a) Let $U\subset(\LL_{\ge 0})^2$ be open. 
       If $U$ intersects the diagonal in a stationary subset, then
       $U\supset [\alpha,\omega_1)^2$ for some $\alpha$.
       If $U$ intersects a horizontal (resp. vertical) line in a stationary subset, then
       $U$ contains a horizontal strip $(\alpha,\omega_1)\times(a,b)$ 
       (resp. vertical strip $(a,b)\times(\alpha,\omega_1)$) for some $\alpha\in\omega_1$
       and $a<b$ in $\LL_{\ge 0}$. \\
   (b) A non-Lindel\"of closed subset of $(\LL_{\ge 0})^2$ intersects either the diagonal, an horizontal line or a vertical
       line on a club subset.\\
   (c) If $f:\mathbb{O}\to\LL_{\ge 0}$ is unbounded on any horizontal line, then it is unbounded on each
       horizontal line and on $\Delta$.\\
   (d) If $U\subset\wh{\mathbb{O}}$ is open and intersects $\Delta$ or $\VV$ in a stationary subset, then
       $\wb{U}\supset [\alpha,\omega_1]\times[\alpha,\omega_1)\cap\wh{\mathbb{O}}$ for some $\alpha$.\\
   (e) Given $g:\wh{\mathbb{O}}\to\R$, there is $c\in\R$ such that 
       $g^{-1}(\{c\})\supset [\alpha,\omega_1]\times[\alpha,\omega_1)\cap\wh{\mathbb{O}}$.\\
   (f) If $g:\wh{\mathbb{O}}\to\R$ is such that $g^{-1}((a,b))\cap \VV$ is non-Lindel\"of, then
       $g^{-1}([a,b])$ contains a terminal part of $\Delta$.\\
   (g) A non-Lindel\"of closed subset of $\wh{\mathbb{O}}$ has a club intersection with $\Delta$ or $\VV$.
\end{lemma}
\proof
   Each property is well known and easy, for instance 
   (a) and (b) are essentially proved in \cite[Lemma 3.4]{Nyikos:1984} or 
   \cite[Lemma B.29 p. 194]{Gauldbook}, (d), (e), (f) and (g) in \cite[8L p.158]{GillmanJerisonBook}.
   For (c), notice that (a) (for horizontal lines) 
   implies that the set $\{ y\in\LL_{\ge 0}\,:\,f\text{ is bounded on }[y,\omega_1)\times\{y\}\}$
   is clopen in $\LL_{\ge 0}$, and hence either empty or the whole space. 
   Applying (a) again (for the diagonal) yields the result.
\endproof

\begin{example}\label{ex:LxR}
   $\LL_{\ge 0} \times [0,1]$ and $\LL_{\ge 0} \times \R$ are f-Type I and narrow (in themselves).
\end{example}
\proof[Details]
Both are obviously f-Type I (use the horizontal projection). 
Assume $\{U_\alpha\,:\,\alpha\in\omega_1\}$ is a systematic cover of $\LL_{\ge 0} \times [0,1]$.
Suppose that $\wb{U_\alpha}$ is non-Lindel\"of for some $\alpha$.
We show that this implies that $\LL_{\ge 0} \times [0,1]$ is contained in some $U_\delta$,
which entails its narrowness.
By Lemma \ref{lemma:propoct} (a) and (b), 
$U_{\alpha+1}$ contains a strip $(\beta_{\alpha},\omega_1)\times(a_\alpha,b_\alpha)$
for $\beta_\alpha\in\omega_1$, $a_\alpha,b_\alpha\in[0,1]$.
Since $\wb{U_{\alpha+1}}\supset [\beta_\alpha,\omega_1)\times\{a_\alpha,b_\alpha\}$, 
$U_{\alpha+2}$  contains a strictly wider (but maybe shorter) strip $(\beta_{\alpha+1},\omega_1)\times(a_{\alpha+1},b_{\alpha+1})$,
that is, $\beta_{\alpha+1}\ge \beta_{\alpha}$ and $a_{\alpha+1}< a_{\alpha}<b_{\alpha}<b_{\alpha+1}$.
Proceeding by induction we obtain an $\omega_1$-sequence of wider and wider strips.
But since $[0,1]$ is hereditarily Lindel\"of, the vertical projection of these strips
must cover all of it after at most countably many steps, hence
some $U_\gamma$ contains $(\beta,\omega_1)\times[0,1]$ for some $\gamma,\beta$.
Now, the space not covered by $U_\gamma$ is compact, hence $\LL_{\ge 0} \times [0,1]$ is contained 
in $U_\delta$ for some $\delta>\gamma$. The proof
$\LL_{\ge 0} \times \R$ is the same.
\endproof

The fact that $[0,1]$ and $\R$ are connected is crucial, as $\LL_{\ge 0}\times X$ is never narrow in itself
when $X$ is the disjoint union of two (non-empty) clopen subsets.

\begin{example}
   $\mathbb{O}$ is f-Type I but not f-narrow, $\wh{\mathbb{O}}$ is f-Type I and narrow.
\end{example}
\proof[Details]
   For $\mathbb{O}$, consider the horizontal and vertical projections:
   the former yields f-Type I-ness, the latter the non-f-narrowness.
   For $\wh{\mathbb{O}}$, the vertical projection shows that it is a f-Type I space.
   Let $\{U_\alpha\,:\,\alpha\in\omega_1\}$ be a systematic cover of $\wh{\mathbb{O}}$.
   If some $\wb{U_\alpha}$ is non-Lindel\"of, 
   by Lemma \ref{lemma:propoct} (g) and (d), it intersects either $\Delta$ or $\VV$ on a club set and
   thus $\wb{U_{\alpha+1}}$ contains 
   $[\alpha,\omega_1]\times[\alpha,\omega_1)\cap\wh{\mathbb{O}}$ for some $\alpha$.  
   Since the rest of $\wh{\mathbb{O}}$ is compact,
   $\wh{\mathbb{O}}$ is contained in $U_{\gamma}$ for some $\alpha<\gamma<\omega_1$.
\endproof

Being (f-)narrow is quite sensitive to the ambient space:
a space might be (f-)narrow in some spaces but not in others. 
Obviously,
given any space, a closed copy of it is narrow in any narrow-in-itself space.
But a non-narrow-in-itself space can become narrow in another non-narrow-in-itself space.

\begin{example}
   Any copy of $\LL\times[0,1]$ in $\mathbb{O}$ is narrow in it, although neither $\LL\times[0,1]$ nor $\mathbb{O}$
   are f-narrow in themselves.
\end{example}
\proof[Details]
   By Lemma \ref{lemma:propoct} (a), any such copy must be contained in $(\LL_{\ge 0}\times[0,\alpha])\cap\mathbb{O}$
   for some $\alpha$. But the latter is Type I and narrow in itself (Example \ref{ex:LxR}).
\endproof

Notice however that it is not true for $\LL$ itself, as $\Delta\cup \HH$ is a copy of $\LL$
not f-narrow in $\mathbb{O}$, as the vertical projection
shows.
Our next result gives many properties equivalent to f-narrowness.

\begin{thm}\label{thm:fnarrow}
  $D\subset X$ be closed and f-Type I in $X$. Then the following properties are all equivalent.\\
  (a) $D$ is f-narrow in $X$.\\
  (b) Given a map $f:X\to\LL_{\ge 0}$ and $E_1,E_2\subset D$ two closed non-Lindel\"of subsets, 
      then $f\upharpoonright E_1$ is bounded iff $f\upharpoonright E_2$ is bounded.\\
  (c) There is no map $f:X\to \LL_{\ge 0}$ with $f\upharpoonright{D}$ unbounded and 
      $f^{-1}(\{0\})\cap D$ non-Lindel\"of.\\
  (d) Let $f:X\to\LL_{\ge 0}$ be such that $f\upharpoonright{D}$ is unbounded and
      $\{U_\alpha\,:\,\alpha\in\omega_1\}$ be a
      systematic chain cover of $X$ with $\wb{U_\alpha}\cap D$ Lindel\"of for each $\alpha$.
      Then for each
      $\alpha\in\omega_1$ there is $\gamma(\alpha)\in\omega_1$ such that 
      $f(D -  U_{\gamma(\alpha)})\subset[\alpha,\omega_1)$.\\
  (e) Let $\{U_\alpha\,:\,\alpha\in\omega_1\}$ be a
      systematic chain cover of $X$ with $\wb{U_\alpha}\cap D$ Lindel\"of for each $\alpha$.
      Then for every $f:X\to\LL_{\ge 0}$ with $f\upharpoonright{D}$ unbounded there is a club $C\subset\omega_1$ such that 
      $f(D -  U_{\alpha})\subset[\alpha,\omega_1)$ for each $\alpha\in C$.\\
  (f) Given $f:X\to\LL_{\ge 0}$, then there is $\alpha\in\omega_1$ such that for any $\beta\ge\alpha$ either 
          $f^{-1}([0,\beta])\cap D$ or $f^{-1}([\beta,\omega_1))\cap D$ is (empty or) Lindel\"of.
\end{thm} 
\proof
   If $D$ is Lindel\"of, each point holds trivially. We thus assume throughout the proof that $D$ is non-Lindel\"of. \\
   (a) $\rightarrow$ (b). Assume that $D$ is f-narrow, and $f$, $E_1$, $E_2$ be as in (b).
      We show that $f\upharpoonright E_i$ is bounded iff $f\upharpoonright D$ is bounded, for $i=1,2$.
      If $f\upharpoonright E_i$ is bounded, then $f^{-1}([0,\alpha])\cap D\supset E_i$ for some $\alpha$.
      By f-narrowness, $f\upharpoonright D$ must then be bounded.
      The conversely is immediate.\\
   (b) $\rightarrow$ (a). Let $f:X\to\LL_{\ge 0}$ and $\alpha\in\omega_1$ be given, and set
      $E_1 = D$ and $E_2 = f^{-1}([0,\alpha])\cap D$.
      If $E_2$ is non-Lindel\"of (b) implies that $f\upharpoonright D$ is bounded.
      This
      shows that $D$ is f-narrow in $X$.\\
   (a) $\rightarrow$ (c). Immediate.
   \\
   (c) $\rightarrow$ (a). 
   If $D$ is not f-narrow, then there is a
   $g:X\to\LL_{\ge 0}$ unbounded on $D$ and an $\alpha$ such that $g^{-1}([0,\alpha])\cap D$ is non-Lindel\"of. 
   Let $p:\LL_{\ge 0}\to\LL_{\ge 0}$
   be continuous such that $p([0,\alpha])={0}$ and $p\upharpoonright{[\alpha+1,\omega_1)}=id$, then $f=p\circ g$ yields a contradiction.
   \\
   (b) $\rightarrow$ (d). 
   If there is some $\alpha$ such that $E_1=f^{-1}([0,\alpha])$ contains points in $X-U_{\gamma}$ for arbitrarily high $\gamma$,
   then $E_1$ is non-Lindel\"of. Since $f$ is unbounded on $D=E_2$, this contradicts (b).\\
   (d) $\rightarrow$ (c). Immediate.\\
   (d) $\rightarrow$ (e). By a leapfrog argument and continuity of $f$. \\
   (e) $\rightarrow$ (d). Immediate.\\
   (b) $\rightarrow$ (f). If $f:X\to\LL_{\ge 0}$ is bounded on $D$, then, for some 
         $\alpha$, $f^{-1}([\alpha,\omega_1))\cap D$ is empty.
         If $f\upharpoonright D$ is unbounded, then $f^{-1}([\alpha,\omega_1))\cap D$ is non-Lindel\"of for each $\alpha$.
         By (b) $f^{-1}([0,\alpha])\cap D$ must be Lindel\"of.\\
   (f) $\rightarrow$ (c). Let $f\co X\to\LL_{\ge 0}$ be unbounded on $D$ and let $\alpha$ be as in (f).
         Since $f^{-1}([\alpha,\omega_1))\cap D$ is non-Lindel\"of, $f^{-1}([0,\alpha])$ is Lindel\"of,
         and hence so is $f^{-1}(\{0\})$.
\endproof

Type I and f-Type I spaces coincide in the realm of regular spaces, this is not the case for f-narrow and narrow spaces for which 
a stronger property (e.g. normality) is needed.
\begin{lemma}\label{lemma:normalfnar}
   Let $X$ be a normal space and $D$ be closed in $X$. Then $D$ is narrow in $X$ iff $D$ is f-narrow in $X$.
\end{lemma}
\proof
   Only the reverse implication needs a proof, and it follows directly from 
   Lemma \ref{lemma:Uknormal}.
\endproof

\begin{figure}
  \begin{center}
    \epsfig{figure = 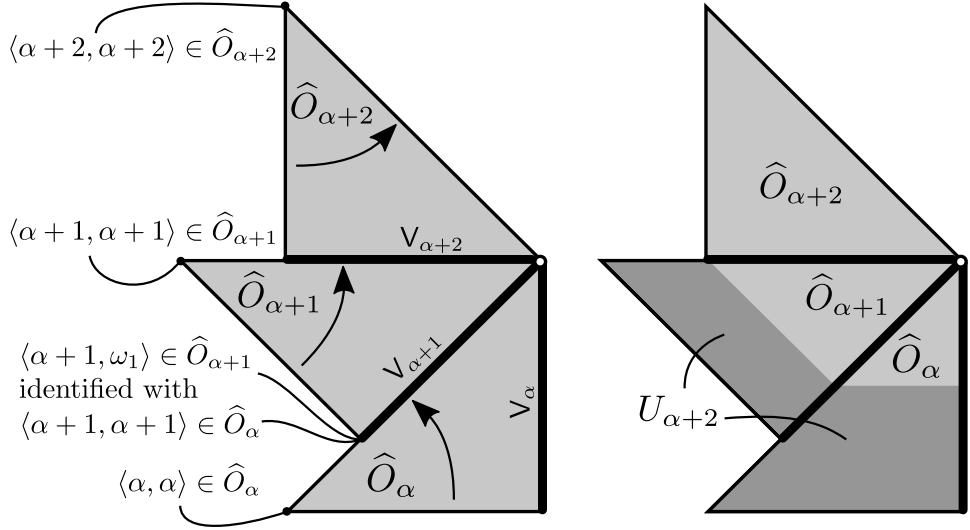, width = .7\textwidth}
    \caption{Piling up copies of $\mathbb{O}_{\ge\alpha}$}
    \label{fig:piling}
  \end{center}
\end{figure}

\begin{example}\label{ex:narnotfunc}
   There is a Type I Tychonoff countably compact locally compact space $Y$ which is f-narrow but 
   not narrow (in itself).
\end{example}
\proof[Details]
   The ideas are similar to that of Example \ref{ex:confuncHaus} (but the construction is simpler).
   As a set, $Y$ is the union of (non disjoint) copies $\wh{O}_\alpha$ of $\wh{\mathbb{O}}_{\ge\alpha}$ for each $\alpha\in\omega_1$.
   Given $\langle x,y\rangle\in\wh{\mathbb{O}}_{\ge\alpha}$, we write $\langle x,y\rangle_\alpha$ for its 
   copy in $\wh{O}_\alpha$, but (somewhat abusively) still write $\VV_\alpha,\Delta_\alpha$ for their respective copies
   in $\wh{O}_\alpha$.
   The topology of each $\wh{O}_\alpha$ outside of $\Delta_\alpha$
   and $\VV_\alpha$ is that of $\wh{\mathbb{O}}_{\ge\alpha}$.
   We glue $\wh{O}_{\alpha+1}$ ``on top'' of $\wh{O}_\alpha$
   by identifying $\langle x,x\rangle_\alpha\in\Delta_\alpha\subset\wh{O}_\alpha$
   with $x\ge\alpha+1$ with $\langle\omega_1,x\rangle_{\alpha+1}\in\VV_{\alpha+1}\subset\wh{O}_{\alpha+1}$,
   as seen on Figure \ref{fig:piling} (left).
   For $\alpha$ limit, a neighborhood of $\langle\omega_1,x\rangle_\alpha\in \VV_\alpha \subset\wh{O}_\alpha$
   is given by fixing $z<x<y$ in $\LL_{\ge 0}$ and $\gamma<\alpha$ in $\omega_1$ and taking
   the union for $\gamma<\beta<\alpha$ of the horizontal bands between height $z,y$ in $\wh{O}_\beta$, 
   together with an open set in $\wh{O}_\alpha$ whose intersection
   with $\VV_\alpha$ is the open segment between $\langle\omega_1,z\rangle_\alpha$ and $\langle\omega_1,y\rangle_\alpha$
   (or the segment $\{\omega_1\}\times[\alpha, z)$ if $x=\langle\alpha,\alpha\rangle\in\wh{O}_\alpha$).
   It should be clear $\langle\alpha,\alpha\rangle_\alpha$
   is the limit point of the sequence $\langle\alpha_n,\alpha_n\rangle_{\alpha_n}$
   for any sequence of ordinals $\alpha_n\nearrow\alpha$. Actually, it is not hard to see that 
   $Y$ is countably compact and that the following claim holds.
   \begin{claim}\label{claim:3.10.1}
      If $z_n\in Y$ ($n\in\omega$) is a sequence such that $z_n = \langle x_n,y_n\rangle_{\alpha_n}$ and
      $\alpha_{n+1} > \max(y_n,\alpha_n)$,
      then this sequence converges to $\langle \alpha,\alpha\rangle_{\alpha}$
      for $\alpha = \sup_{n\in\omega}\alpha_n$.
   \end{claim}
   Define $\mathsf{W}$ as $\{\langle \alpha,\alpha\rangle_{\alpha}\,:\,\alpha\in\omega_1\}$.
   By the previous claim, $\mathsf{W}$ is a copy of $\omega_1$ in $Y$.
   Let $E_{\alpha,\beta}$ be 
   $\{\langle x,y\rangle_\alpha\in\wh{O}_\alpha\,:\,y\ge\beta\}$, that is, 
   the points in $\wh{O}_\alpha$ with vertical coordinate $\ge\beta$.
   \begin{claim}\label{claim:3.10.2}
      Let $C\subset Y$. 
      (i)
      If there is some $\beta$ such that 
      $C\cap E_{\gamma,\beta}\not=\varnothing$ for a subset of $\gamma$ which is unbounded below some 
      limit $\alpha<\omega_1$, then $\wb{C}\cap\VV_\alpha\cap E_{\alpha,\beta}\not=\varnothing$. 
      (ii) If the set of $\gamma<\omega_1$ such that $C\cap \wh{O}_{\gamma}\not=\varnothing$
      is unbounded in $\omega_1$,
      then $\wb{C}\cap\mathsf{W}$ is club.
   \end{claim}
   \proof (i) Fix $\gamma_n\nearrow\alpha$ and $\langle x_n,y_n\rangle_{\gamma_n} \in E_{\gamma_n,\beta}\cap C$.
          Then take a converging subsequence to obtain a point in $\VV_\alpha$.
          (ii) For $\mathsf{W}$, fix $\gamma_0\in\omega_1$,
          and by induction choose $\langle x_n,y_n\rangle_{\alpha_n}\cap C$
          such that $\alpha_{n+1} > \max(y_n,\alpha_n)$.
          Then apply Claim \ref{claim:3.10.1}.
   \endproof
   It should be clear that $Y$ is locally compact (and thus Tychonoff).
   Moreover, $Y$ is Type I (functionally, thanks to Lemma \ref{lemma:TypeIregular}) with canonical cover 
   $$ U_\alpha = \cup_{\beta<\alpha} \{\langle x,y\rangle_\beta\in \wh{O}_{\beta}\,:\, y<\alpha\}$$
   (see Figure \ref{fig:piling}, right).
   Actually, $\wb{U_\alpha}$ is compact for each $\alpha$.
   The systematic cover
   $$ \{(\cup_{\beta<\alpha} \wh{O}_{\beta}) - \VV_\alpha\,:\,\alpha\in\omega_1\}  $$
   shows that $Y$ is not narrow it itself.
   However, $Y$ is functionally narrow in itself, as we shall now see.
   \begin{claim}\label{claim3}
      (i) For any open set $U$ whose intersection with $\mathsf{V}_{\alpha}$ is stationary for some $\alpha<\omega_1$ 
      there are $\gamma<\alpha$ and $\delta<\omega_1$ such that $\wb{U}$ contains
      $\cup_{\gamma\le\beta<\alpha} E_{\beta,\delta}$.
      (ii) Any open set whose intersection with $\mathsf{W}$ is stationary 
      (that is, it contains $\langle\alpha,\alpha\rangle_\alpha$ for a stationary set of $\alpha$)
      contains $\cup_{\beta>\gamma}\wh{O}_{\beta}$
      for some $\gamma$.
   \end{claim}
   \proof
      We show (ii) first.
      Let $U$ be open and intersect $\mathsf{W}$ on a stationary set and assume that the claim does not hold.
      Then the complement of $U$ has a club intersection with $\mathsf{W}$ by Claim \ref{claim:3.10.2}
      and thus intersects $U$, a contradiction.
      (i)
      If $\alpha=\gamma+1$, then $\wb{U}$ contains $E_{\gamma,\beta}$ for some $\beta$
      by Lemma \ref{lemma:propoct} (d). 
      If $\alpha$ is limit, it follows from the previous claim.
   \endproof 
   \begin{claim}
      If $g\co Y\to\R_{\ge 0}$ is such that $g^{-1}([0,a])$ is non-Lindel\"of for some $a\ge0$,
      then $g^{-1}([0,b])\cap\mathsf{W}$ is club for any $ b>a$
      and contains each $\wh{O}_\alpha$ above some $\gamma$.
   \end{claim} 
   \proof
      Let $b>a$ be given.
      If $g^{-1}([0,b))\cap\mathsf{W}$ is stationary, we are over by Claim \ref{claim3}.
      If not, by Claim \ref{claim:3.10.2} there is $\gamma<\omega_1$ such that $g^{-1}([0,b))\cap\wh{O}_\alpha=\varnothing$ 
      for each $\alpha\ge\gamma$.
      Hence $g^{-1}([0,a])\cap\wh{O}_\alpha$ is non-Lindel\"of for some $\alpha<\gamma$.
      By Lemma \ref{lemma:propoct} (d) \& (f), $g^{-1}([0,c])$ contains a terminal part of $\mathsf{V}_\alpha$ for any $c>a$.
      Let $e\co \gamma\to\R$ be a monotone embedding of the ordinal $\gamma+1$ such that $e(\alpha) = a$
      and $a<e(\gamma)<b$.
      Applying inductively Lemma \ref{lemma:propoct} (d) \& (f), we see that 
      $g^{-1}([0,e(\beta)))$ contains a terminal part of $\mathsf{V}_\beta$ for each
      $\beta$ between $\alpha$ and $\gamma$ (included). Hence, $g^{-1}([0,b))$ must contain a terminal part
      of $\mathsf{V}_\gamma$, contradicting $g^{-1}([0,b))\cap\wh{O}_\gamma=\varnothing$.
   \endproof
   Now let $f\co Y\to\LL_{\ge 0}$ be given with $f^{-1}([0,\alpha])$ non-Lindel\"of.
   Set $g(x) = \min\{f(x), \alpha+\omega\}$, then 
   $g$ has range $\subset[0,\alpha+\omega+1)\simeq\R_{\ge 0}$, and
   by the previous claim $f^{-1}([0,\alpha+1])=g^{-1}([0,\alpha+1])$
   contains each $\wh{O}_\alpha$ above some $\gamma$.
   By Claim \ref{claim3} applied recursively, there is $n\in\omega$ and $\beta<\omega_1$ such that $f^{-1}([0,\alpha+n])$
   contains $\cup_{\alpha < \gamma} E_{\alpha,\beta} \bigcup\cup_{\gamma\le\alpha<\omega_1} \wh{O}_\alpha$.
   But the complement of this set is compact, hence $f^{-1}([0,\delta])$ contains all of $Y$ for some $\delta\in\omega_1$
   bigger than $\alpha+n$,
   that is $f$ is bounded. This shows that $Y$ is f-narrow.
\endproof

\begin{q}
   Is there a first countable space as in Example \ref{ex:narnotfunc}~?
\end{q}

%%%%%%%%%%%%%%%%%%%%%%%%%%%%%%%%%%%%%%%%%%%%%%%%%%%%%%%%%%%%%%%%%%%%%%%%%%%%%%%%%%%%%%%%%%%%%%%%%%%%%%%%%%%%%%%%%%%%%%%%%%%%%%%%%%%

\section{Type I spaces without (f-)narrow non-Lindel\"of subspaces}\label{sec:without}

There are trivial examples of f-Type I spaces without non-Lindel\"of narrow subspaces.

\begin{example}
   A discrete space of cardinality $\aleph_1$ is f-Type I and does not contain any non-Lindel\"of subspace f-narrow in it.
\end{example}

It is more difficult to find examples if one imposes further properties on the space such as 
countable compactness or $\omega_1$-compactness.
If one adds to the mix some local demands (first countability or weakanings of it), 
the existence of such spaces ends up depending on the axioms of set theory.
Let us first show that there is a countably compact f-Type I space without non-Lindel\"of subspaces f-narrow in it.
This example is not countably tight, as none can be exhibited in {\bf ZFC} alone since 
the proper forcing axiom {\bf PFA} prevents their existence
(see Theorem \ref{thm:PFA1} below). Recall that $\beta\omega$ is the \v Cech-Stone compactification of the integers.

\begin{example}
   There is an $\omega$-bounded f-Type I subspace $Q$ of $\omega^\ast=\beta\omega-\omega$ 
   containing no closed non-compact subspace f-narrow in $Q$.
\end{example}
Recall that $\omega$-boundedness and countable compactness are equivalent for Type I spaces.
\proof[Details]
  Let $\mathcal{A}=\{A_{\alpha}\,:\,\alpha\in\omega_1\}$ be a $\subset^\ast$ strictly increasing $\omega_1$-sequence of
  infinite subsets of $\omega$, that is: $A_\alpha-A_\beta$ is finite and $A_\beta-A_\alpha$ infinite when
  $\alpha < \beta$. (It is a classical fact that $\mathcal{A}$ can be defined by induction.)
  Set $U_{\alpha} = \wb{A_{\alpha}}-\omega$, where the closure is taken in $\beta\omega$.
  Then each $U_\alpha$ is clopen in $\omega^\ast$
  and $U_\beta$ contains each $U_\alpha$ whenever $\alpha<\beta$ (see e.g. \cite[6Q {\&} 6S]{GillmanJerisonBook}).
  Hence, $Q=\cup_{\alpha<\omega_1}U_\alpha$ is a Tychonoff Type I-in-itself (and hence f-Type-I-in-itself) subspace of 
  $\omega^\ast$. Notice that $Q$ is open in $\omega^\ast$.
  A canonical cover of $Q$ is given by replacing $U_\alpha$ by  $\cup_{\beta<\alpha}U_\beta$ at limit levels.
  Since each $U_{\alpha+1}$ is compact, $Q$ is $\omega$-bounded. Let $C\subset Q$ be closed and non-Lindel\"of.
  We show that $C$ is not f-narrow in $Q$.
  Up to taking a subcover of the canonical cover, we may assume that 
  $C\cap(U_{\alpha+1}-\wb{U_\alpha})$ contains at least two points.
  We shall define closed subsets 
  $D_\alpha,E_\alpha\subset C\cap \wb{U_\alpha}$ and $f_\alpha:\beta\omega\to[0,\alpha]\subset\LL_{\ge 0}$
  such that the following holds:
  \\
  (1) $f_\alpha(D_\alpha)=\{0\}$, $f_\alpha(E_\alpha)=\{0,1,2,\dots,\alpha\}\subset\LL_{\ge 0}$),
  \\
  (2) $f_\alpha\upharpoonright{\wb{U_\gamma}}=f_\gamma\upharpoonright{\wb{U_\gamma}}$ whenever $\gamma\le\alpha$,
  \\
  (3) $D_\alpha\cap \wb{U_\gamma}= D_\gamma$ and $E_\alpha\cap \wb{U_\gamma}= E_\gamma$ whenever $\gamma\le\alpha$.\\
  Once these are defined, set $D=\cup_{\alpha<\omega_1}D_\alpha$, $E=\cup_{\alpha<\omega_1}E_\alpha$.
  Then the map $f:Q\to\LL_{\ge 0}$ defined by $f(x)=f_{\alpha}(x)$ for $\alpha$ such that $x\in U_\alpha$
  is continuous, unbounded on $E$ and constant on $0$ on $D$. This shows that 
  $C$ is not f-narrow in $X$ by Theorem \ref{thm:fnarrow} (c).
  \\
  So, start with $E_0=F_0=\{e_0\}\subset U_0$, $f_0\co U_0\to\LL_{\ge 0}$ be constant on $0$.
  Let $\alpha=\gamma+1$.
  Take $d,e\in C\cap(U_{\gamma+1}-\wb{U_\gamma})$. The map $g:U_\gamma\cup \{d,e\}\to[0,\gamma+1]\subset\LL_{\ge 0}$ 
  defined as 
  $$
      g(x)=\left\{
          \begin{array}{cl}
            f_\gamma(x) &\text{if }x\in \wb{U_\gamma} \\
            0 & \text{if }x=d \\ 
            \gamma+1 &  \text{if }x=e
          \end{array}\right.
  $$
  is continuous. 
  Since $\wb{U_\gamma}\cup \{d,e\}$ is closed and $\beta\omega$ normal,
  there is some $f_{\gamma+1}:\beta\omega\to[0,\gamma+1]$ extending $g$. Set $D_{\gamma+1}=D_\gamma\cup\{d\}$,
  $E_{\gamma+1}=E_\gamma\cup\{e\}$, by construction the above conditions are fulfilled.
  Let now $\alpha$ be limit. The map $g:U_\alpha\to[0,\alpha]$ defined by 
  $g\upharpoonright{U_\gamma}=f_\gamma\upharpoonright{U_\gamma}$ 
  for $\gamma<\alpha$ is continuous. 
  As well known (see e.g. \cite[Theorem 1.5.2]{vanMill:1984}),
  if $M\subset\beta\omega$ is Lindel\"of, then
  any continuous function $M\to[0,1]$ can be continuously extended 
  to a function $\beta\omega\to[0,1]$.
  Hence, there is 
  $f_\alpha:\beta\omega\to[0,\alpha]$ extending $g$. Set 
  \begin{align*}
      D_\alpha&=(\cup_{\gamma<\alpha}D_\gamma)\cup\bigl(f_\alpha^{-1}(\{0\})\cap C \cap (\wb{U_\alpha}-U_\alpha)\bigr), \\
      E_\alpha&=(\cup_{\gamma<\alpha}E_\gamma)\cup\bigl(f_\alpha^{-1}(\{\alpha\})\cap C \cap (\wb{U_\alpha}-U_\alpha)\bigr).
  \end{align*}
  By construction $D_\alpha$ and $E_\alpha$ are closed in $X$ and satisfy the above conditions (1)--(3).
  This finishes the proof.
\endproof

We now show that there are (very classical) consistent examples of $\omega_1$-compact first countable
Type I spaces without any closed non-Lindel\"of subspace f-narrow in it, namely: {\em Suslin trees}.
Recall that a Suslin tree is an tree of height $\omega_1$ whose chains and antichains are at most countable.
They cannot be shown to exist by {\bf ZFC} alone but pop up in many models of set theory.
In what follows, a
(set theoretic) tree $T$ is always given the order (also called interval) topology,
that is, a basis for the open sets is given by the ``intervals'' $I_{x,y}=\{ z\in T\,:\,x<z<y\}$.
Any tree is locally compact and a tree of height $\le\omega_1$ is first countable in this topology.
We shall often use the fact that antichains are closed discrete in trees, and that a 
closed discrete subset of a Hausdorff tree is a countable union of antichains 
(see e.g. Theorem 4.11 in \cite{Nyikos:trees}). A Suslin tree is $\omega_1$-compact.
Recall also that a tree is Hausdorff
iff two members at a limit level are equal whenever they have the same predecessors.
The next lemma shows a little bit more than what we really need.

\begin{lemma}\label{lemma:treenotnarrow}
   A Hausdorff tree $T$ of height $\omega_1$ 
   without an uncountable branch does not contain any closed
   non-Lindel\"of subspace which is both Type I in $T$ and f-narrow in $T$.
\end{lemma}
\proof
   Let $D\subset T$ be closed. If $D$ has points as high as one wants in the tree,
   since $T$ has no uncountable branch, there are $x_0,x_1\in T$ at same level $\alpha$ such that $\{y\in D\,:\, y>x_i\}$
   is uncountable for both $i=0,1$.  
   Set $f\co T\to\LL_{\ge 0}$ to be $0$ until level $\alpha$ and above $x_0$, and equal to the level function elsewhere.
   Then $f$ shows that $D$ is not functionally narrow.\\
   If $D$ is bounded in the tree, then $D$ is either countable and hence Lindel\"of,
   or its intersection with some level of $T$ is uncountable.
   Let $N$ be the set of $\alpha$ such that 
   $A_\alpha = \{ u\in\text{Lev}_\alpha(T)\,:\,\exists d\in D\text{ with } u \le d\}$
   is uncountable. Notice that  when $\alpha\in N$, $A_\alpha$
   is an uncountable closed discrete subspace of $T$.
   If $N$ contains a successor ordinal $\alpha$ (e.g. if it contains more than one element),
   split $A_\alpha$ in two disjoint uncountable antichains $B,C$.
   Let $f\co B\to\omega_1\subset\LL_{\ge 0}$ be onto. Extend it to all of $T$ by taking value $0$ outside of $B$.
   Then $f$ is continuous and shows that $D$ is not f-narrow.
   \\
   If $N$ contains only one element $\alpha$ at a limit level, then $A_\alpha\subset D$ and
   all levels of $T_0 = \{ t\in T\,:\, \exists d\in A_\alpha\text{ with }t< d\}$
   are countable.
   Thus $T_0$ itself is countable,  
   and for any systematic chain cover $\{U_\beta\,:\,\beta\in\omega_1\}$ of $T$,
   $A_\alpha\subset \wb{U_\beta}$ for (the first) $\beta$ such that $T_0\subset U_\beta$,
   hence $D\cap\wb{U_\beta}$ is non-Lindel\"of and $D$ is not Type I in $T$.
\endproof

Recall that an {\em $\omega_1$-tree} is a tree of height $\omega_1$ whose levels are countable.

\begin{lemma}
   Let $D$ be an uncountable subset of a tree $T$ of height $\omega_1$.
   Then either $D^\downarrow = \{t\in T\,:\,\exists d\in D\text{ with }t\le d\}$ is an $\omega_1$-tree,
   or $D$ contains an uncountable antichain.
\end{lemma}
\proof
   If $D^\downarrow$ is not an $\omega_1$-tree, consider its smallest uncountable level.
   Pick one point in $D$ above each member in this level. This yields an uncountable antichain.
\endproof

\begin{lemma}
   Any closed subspace of an $\omega_1$-tree $T$ is f-Type I in $T$.  
\end{lemma}
\proof The level function works.
\endproof

The above lemmas yield immediately our claimed example.

\begin{example}\label{ex:Suslin}
   A Suslin tree is an $\omega_1$-compact locally compact first countable f-Type I space containing no f-narrow non-Lindel\"of subspace.
\end{example}

\begin{q}\label{q:nonarrow}
   Consistently, is there an $\omega$-bounded first countable (f-)Type I space containing no closed
   non-Lindel\"of subspace (f-)narrow in it~?
\end{q}

This question dates back to 2006, see \cite[Problem 6.3]{mesziguessurf}. 
Peter Nyikos, in 2006 also, claimed to have found an $\omega$-bounded surface (i.e. $2$-manifold) 
containing no closed
non-Lindel\"of subspace f-narrow in it under 
$\diamondsuit^+$,
but he seems to have unfortunately been sidetracked by many other problems 
before completing a first draft of his construction\footnote{Another proof of the high level
of involvment of experts in the present subject.}, which we were not able to emulate.
Question \ref{q:nonarrow} asks for a consistent example because of the following theorem.

\begin{thm}[{\cite{Balogh:1989, Eisworth:2002}, in effect}]{\bf (PFA)} \label{thm:PFA1}
    If $X$ is countably compact, countably tight and Type I in some space $Y$, 
    then $X$ contains a closed copy of $\omega_1$ and thus a subspace narrow in it.              
\end{thm}
\proof
  Given a canonical cover $\mathcal{U}$
  of $X$ in $Y$, its skeleton $\text{Sk}(\mathcal{U},X)$ is a perfect preimage of $\omega_1$.
  (Send the $\alpha$-th bone to $\alpha$, this gives a continuous closed map with compact fibers.)
  Under {\bf PFA}, it contains a closed copy of $\omega_1$.
  This was shown by Balogh in \cite{Balogh:1989} for spaces of character $\le\aleph_1$ and 
  Eisworth proved it
  for countably tight spaces in \cite{Eisworth:2002}. 
\endproof
Notice that for first countable spaces, 
the conclusion holds in
a model of {\bf ZFC + CH} without inaccessible cardinals by a result of Eisworth and
Nyikos \cite{EisworthNyikos:2005}.
The paternity (through an oral account) of the next result should be attributed to P. Nyikos alone.

\begin{thm}{\bf (PFA)} 
   If $X$ is $\omega_1$-compact, countably tight, locally compact and Type I in some space $Y$, 
   then $X$ 
   contains a closed copy of $\omega_1$ and thus a subspace narrow in $X$ (and in $Y$).              
\end{thm}
Example \ref{ex:Suslin} shows that this result does not hold in {\bf ZFC}.
\proof
   If $X$ is Lindel\"of, there is nothing to prove, hence we assume that it is not.
   The first Trichotomy Theorem of Eisworth and Nyikos \cite{EisworthNyikos}
   shows that under {\bf PFA}
   a locally compact space is either a countable union of
   $\omega$-bounded subspaces, or has a closed uncountable discrete space, or has a Lindel\"of
   subset with a non-Lindel\"of closure. The latter two are impossible if $X$ is Type I in $Y$ and
   $\omega_1$-compact. Thus $X$ is a countable union of $\omega$-bounded subspaces, one of which must
   be non-Lindel\"of, and we evoke Theorem \ref{thm:PFA1} to conclude.
   (Note that a copy of $\omega_1$ is closed in a countably tight Type I subspace of $Y$.)
\endproof

We note in passing that there are models of set-theory in which there are $\omega$-bounded
first countable spaces ($2$-manifolds, even)
containing no copies of $\omega_1$ but which are however narrow in themselves, for instance \cite[Example 6.9]{Nyikos:1984}
(built with the help of $\diamondsuit$). 
This brings the next question.

\begin{q}
   Is there a model of set theory where first countable $\omega$-bounded Type I spaces neccessarily contain a non-Lindel\"of
   closed subspace (f-)narrow in it, but not necessarily a copy of $\omega_1$~?
\end{q}

Another question, inspired by Example \ref{ex:narnotfunc}.

\begin{q}
   Is there a space which contains a non-Lindel\"of closed subspace f-narrow in it but 
   no closed non-Lindel\"of subspace narrow in it~? Is there an f-Type I example~?
\end{q}

Notice that each subspace $\VV_\alpha$ (and $\mathsf{W}$) is narrow in the space of Example \ref{ex:narnotfunc},
and that any subset is f-narrow in Example \ref{ex:conHaus} 
(since functions with range $\subset\LL_{\ge 0}$ are constant), but this space also contains a copy of $\omega_1$.

%%%%%%%%%%%%%%%%%%%%%%%%%%%%%%%%%%%%%%%%%%%%%%%%%%%%%%%%%%%%%%%%%%%%%%%%%%%%%%%%%%%%%%%%%%%%%%%%%%%%%%%%%%%%%%%%%%%%%%%%%%%%%%%%%%%

\section{Discrete narrow and f-narrow non-Lindel\"of subspaces}\label{sec:discrete}

Recall that a space $X$ is {\em (strongly) collectionwise Hausdorff} ({\em (s)cwH} for short) 
iff any closed discrete subspace $\{d_\alpha\,:\,\alpha\in\kappa\}$
can be expanded to a disjoint (resp. discrete) collection of open sets.
That is, there are open $O_\alpha\ni d_\alpha$ such that $\{O_\alpha\,:\,\alpha\in\kappa\}$
is a disjoint (resp. discrete) collection in $X$.

There are simple examples of spaces containing discrete narrow non-Lindel\"of subspaces,
for instance $\LL_{\ge 0}\times[0,1] - \omega_1\times\{0\}$ where $\omega_1$ is seen as a subspace of $\LL_{\ge 0}$, 
see Example \ref{ex:cwHnarrow} for details. The main goal of this section is to exhibit 
more elaborate examples whose narrow-in-the-space
non-Lindel\"of subspaces are all discrete outside of a Lindel\"of subset.
We first show how scwH-ness does not allow such behaviours and how f-narrow subspaces interact with 
skeletons of canonical covers. At the same time, we investigate how much information about (f-)narrow subspaces does
the skeleton possess.

\begin{lemma}\label{lemma:7.2}
    Let $D$ be a closed subspace of $Y$ and $\{V_\alpha\,:\,\alpha\in\omega_1\}$ be a good systematic cover of $Y$
    such that $\wb{V_\alpha}$ is Tychonoff 
    and $D\not\subset V_\alpha$ for each $\alpha$. Let $U\supset D$ be open. 
    If $D$ is f-narrow in $Y$, 
    then $\wb{U}\cap(\wb{V_\alpha}-V_\alpha)\not=\varnothing$ for a stationary set of $\alpha$. 
\end{lemma}
\proof 
    By contradiction,
    suppose that $\wb{U}\cap(\wb{V_\alpha}-V_\alpha)=\varnothing$ 
    for a club set of $\alpha$. Since $D\not\subset V_\alpha$,
    up to taking a subscover, we may assume that
    that $\wb{U}\cap(\wb{V_\alpha}-V_\alpha)=\varnothing$ and
    $D\cap U_\alpha\not=\varnothing$ for each $\alpha$, where 
    $U_\alpha = U\cap(V_{\alpha+1}-\wb{V_\alpha})$.
    If $y\in\wb{U}$, 
    take $\alpha$ minimal such that $y\in V_\alpha$, then since the cover is good $\alpha$ must be a successor ordinal $\beta+1$.
    Then $y\in V_{\beta+1}-\wb{V_\beta}$, which intersects at most
    one member of the family $\mathcal{U}=\{U_\alpha\,:\,\alpha\in\omega_1\}$, which is thus 
    discrete.
    Picking a point $x_\alpha$ in each $D\cap U_\alpha$ yields
    a club discrete subset of $D$ which we may partition in two disjoints uncountable club discrete subsets 
    $D_i=\{x_\alpha\,:\,\alpha\in E_i\}$, $i=0,1$.
    Since $\wb{V_\alpha}$ is Tychonoff and contains $U_\alpha$, 
    we may define $f_\alpha\co \wb{V_\alpha}\to[0,\alpha]\subset\LL_{\ge 0}$
    which takes value $\alpha$ on $x_\alpha$ and $0$ outside of $U_\alpha$.
    By discreteness of $\mathcal{U}$, the function $f\co Y\to\LL_{\ge 0}$ equal to $f_\alpha$ inside of $U_\alpha$
    for $\alpha\in E_1$
    and $0$ outside of of $\cup_{\alpha\in E_1} U_\alpha$ is continuous.
    But $f^{-1}(\{0\})\cap D$ is non-Lindel\"of since it contains $D_0$, hence
    $D$ is not f-narrow in $Y$.
\endproof

\begin{cor}
   If $X$ is a regular Type I space and $\mathcal{V}$ is a canonical cover of $X$, then the closure of any 
   open set containing a closed
   non-Lindel\"of subset of $X$ which is f-narrow in
   $X$ intersects the members of $\BO(\mathcal{V})$ on a stationary set of indices.
\end{cor}
\proof By Lemma \ref{lemma:TypeIregular} $X$ is Tychonoff and the result follows by the previous lemma.
\endproof

The following lemma has a very similar proof. Recall that a space is {\em ccc} iff
any disjoint collection of open sets is at most countable.
\begin{lemma}\label{lemma:cwHdiscrete}
   Let $D$ be closed discrete in the cwH Type I space $X$ with canonical cover $\mathcal{V}$.
   If members of $\mathcal{V}$ are ccc,
   then $D$ misses the members of $\BO(\mathcal{V})$ on a club set of indices.
\end{lemma}
\proof
   Since $X$ is cwH, there is a disjoint family $\mathcal{O}$ of open sets containing each a point of $D$.
   If $D$ intersects the members of $\BO(\mathcal{V})$ on a stationary set of indices,
   by Fodor's lemma there is some member $V_\alpha\in\mathcal{V}$ intersected by
   uncountably many members of $\mathcal{O}$. This is impossible since $V_\alpha$ is ccc.
\endproof

\begin{lemma}\label{lemma:7.3}
   If $X$ is Tychonoff and scwH then $X$ contains no uncountable closed discrete subspace which is
   f-narrow in $X$.
\end{lemma}
\proof 
   Given an uncountable closed discrete $D$, expand it to a discrete collection of open sets $\mathcal{O}$.
   Partition $\mathcal{O}$ in two uncountable disjoint families $\mathcal{O}_0, \mathcal{O}_1$.
   Then as in Lemma \ref{lemma:7.2} define a fonction $f\co X\to\LL_{\ge 0}$ which 
   is $0$ on the complement of $\cup\mathcal{O}_1$ 
   and takes higher and higher values in each member of $\mathcal{O}_1$.
   Discreteness ensures that the function is continuous.
\endproof

``Strongly'' cannot be omitted in Lemma \ref{lemma:7.3}, as cwH is not enough to rule 
out closed discrete subspaces narrow in the whole space. 

\begin{example}\label{ex:cwHnarrow}
   $Y=\LL_{\ge 0}\times\R-\omega_1\times\Q$ is cwH, f-Type I and narrow in itself.
   Any open non-Lindel\"of subset contains
   an uncountable closed discrete subset that is narrow in $Y$.
\end{example}
Of course, $\omega_1$ is seen here as a subspace of $\LL_{\ge 0}$.
\proof[Details]
   The projection on the first factor ensures that $Y$ is functionally Type I. 
   Let us show that $Y$ is narrow in itself. The proof is very similar to that of Example \ref{ex:LxR}.
   \begin{claim} Let $U,V\subset Y$ be open such that $U$ intersects some horizontal 
                 $\LL_{\ge 0}\times\{y\}\cap Y$ 
                 unboundedly and $V\supset\wb{U}$.
                 Then $V\supset [x,\omega_1)\times\bigl((a,b)-\Q\bigr)$ for some $x\in\LL_{\ge 0}$ and
                 $a<y<b$ in $\R$.
   \end{claim}
   \proof
   If $y\in\R-\Q$, this is a consequence of Lemma \ref{lemma:propoct} (a) (and subspace topology)
   since $\LL_{\ge 0}\times\{y\}\cap Y$ is a copy of $\LL_{\ge 0}$.
   If $y\in\Q$, take a family $\{x_\alpha\in\LL_{\ge 0}\,:\,\alpha\in\omega_1\}$
   such that $x_\alpha>\alpha$ and
   $\langle x_\alpha, y\rangle\in U$.
   Then $U\supset \{x_\alpha\}\times(a_\alpha,b_\alpha)$ for some $a_\alpha,b_\alpha\in\Q$,
   $a_\alpha<y<b_\alpha$, and some interval $(a,b)$ must appear uncountably many times.
   It follows that $\wb{U}$ intersects $C\times[a,b] \cap Y$ for some club $C\subset\LL_{\ge 0}$
   and $a<y<b$ in $\R$. 
   Hence for each $y\in(a,b)-\Q$, $V$ contains $[x_y,\omega_1)\times(a_y,b_y)\cap Y$.
   By Lindel\"ofness, there is a countable family of $y$ such that the union
   of $(a_y,b_y)$ cover $(a,b)-\Q$, taking the supremum of the $x_y$ yields the claim.
   \endproof   
   \noindent Suppose that $\{U_\alpha\,:\,\alpha\in\omega_1\}$ 
   is a systematic cover of $Y$ with $\wb{U_0}$ non-Lindel\"of, hence $U_0$ is horizontally unbounded.
   Since $\R$ is second countable there is 
   some interval $(a,b)$ such that $U_0$ contains $\{x\}\times(a,b)$ for unboundedly many $x$.
   We see by induction and the above claim that 
   $\wb{U_\beta}$ contains $[x_\beta,\omega_1)\times[a_\beta,b_\beta]\cap Y$
   with $x_\beta\ge x_\alpha$ and 
   $a_\beta < a_\alpha < b_\alpha < b_\beta$ whenever $\alpha<\beta$.   
   By secound countability again, there is some $\beta<\omega_1$ such that
   $\cup_{\gamma<\beta}(a_\gamma,b_\gamma) = \R$.
   Hence $\wb{U_\gamma}\supset [x_\gamma,\omega_1)\times\R\cap Y$.
   Since the rest of the space is Lindel\"of, it is contained in some $U_\delta$ for $\gamma<\delta<\omega_1$.
   It follows that $Y$ is narrow in itself.
   Also,
   if $U$ is open and non-Lindel\"of, there is  
   $q\in\Q$ such that $U\cap (\LL_{\ge 0}-\omega_1)\times\{q\}$
   is unbounded. It is then easy to find an uncountable closed discrete subset of 
   $U\cap (\LL_{\ge 0}-\omega_1)\times\{q\}$.
   To finish, we are left only with the following claim to prove.
   \begin{claim}
     $Y$ is cwH.
   \end{claim}
   \proof By \cite[Lemma B.29 p. 194]{Gauldbook} (or arguing as above) 
   a closed
   subset $C$ of $\LL_{\ge 0}\times(\R-\Q)$ which has unbounded projection on the first factor
   must intersect $\LL_{\ge 0}\times\{c\}$ on a club set for some $c\in\R-\Q$.
   If $C$ is discrete, its intersection with $\LL_{\ge 0}\times(\R-\Q)$ must then be 
   horizontally bounded, say by $\gamma$.
   Since $Y\cap [0,\gamma+2)\times\R$ is metrizable, 
   it is scwH and we may separate the points of $C\cap [0,\gamma]\times\R$
   by open sets inside of it. The rest of $C$ lies in 
   $\displaystyle\cup_{\gamma<\alpha<\omega_1}(\alpha,\alpha+1) \times \R$,
   which is a disjoint (non-discrete) union of open metrizable subspaces of $Y$.
   We may thus separate the rest of $C$ in each piece by open sets.
   \endproof
   \noindent 
   This finishes the proof
\endproof

To obtain more, that is, spaces $X$ where all narrow-in-$X$ non-Lindel\"of
subspaces are discrete after removal of
a Lindel\"of subset (and some other variations), we introduce 
a general construction procedure to be applied later in concrete cases.

\begin{defi}
   Let $X$ be f-Type I with $s\co X\to\LL_{\ge 0}$ a slicer of $X$.
   We set:
   \begin{align*}
      \text{Down}(X,s)& = \{\langle x,y\rangle\in X\times\LL_{\ge 0}\,:\, y\le s(x)\}\\
      \Delta(X,s) &= \{\langle x,y\rangle\in X\times\LL_{\ge 0}\,:\, y = s(x)\}
   \end{align*}
\end{defi}
   $\Delta(X,s)$ is homeomorphic to $X$ (see below), and $\text{Down}(X,s)$ is a closed subspace of $X\times\LL_{\ge 0}$.
   If $X=\LL_{\ge 0}$ and $s$ is the identity, then $\text{Down}(X,s)$ is the octant $\mathbb{O}$
   and $\Delta(X,s)$ its diagonal $\Delta$.

For our later examples, we only need the case where $X$ is an $\omega_1$-tree
and $s$ is the height function, but it is not 
difficult to obtain general properties of $\text{Down}(X,s)$.

\begin{lemma}\label{lemma:down}
   Let $X$ be f-Type I with slicer $s\co X\to\LL_{\ge 0}$, $D\subset \text{Down}(X,s)$ be closed
   and $\pi\co \text{Down}(X,s)\to X$ be the projection
   on the first factor. Then the following hold.\\
   (a) $\pi$ is a closed map and $\pi$ restricted to $\Delta(X,s)$ is a homeomorphism.\\
   (b) $s\circ\pi$ is a slicer of $\text{Down}(X,s)$.\\
   (c) If there is no closed non-Lindel\"of subset f-narrow in $X$, then there is no closed non-Lindel\"of subset f-narrow in
       $\text{Down}(X,s)$.\\
   (d) If $D$ is closed discrete, so is $\pi(D)$.\\
   (e) If $\pi(D)$ contains a closed discrete subset, then so does $D$.\\
   (f) $X$ is countably compact iff $\text{Down}(X,s)$ is countably compact.\\
   (g) $X$ is $\omega_1$-compact iff $\text{Down}(X,s)$ is $\omega_1$-compact.\\
   (h) $X$ is locally compact iff $\text{Down}(X,s)$ is locally compact.
\end{lemma}
\proof\
   \\
   (a) The homeomorphism claim is immediate since $\Delta(X,s)$ is the graph of $s$ in $X\times\LL_{\ge 0}$.
   Set $U_\alpha = s^{-1}([0,\alpha))$ for each $\alpha$.
   Let $C\subset\text{Down}(X,s)$ be closed and $x\in \wb{\pi(C)}\subset X$.
   Let $\alpha$ be such that $U_\alpha\ni x$.
   Then $x$ is in the closure of $\pi\bigl(C\cap (\wb{U_\alpha}\times[0,\alpha])\bigr)$.
   By compactness of $[0,\alpha]$
   the projection $\wb{U_\alpha}\times[0,\alpha]\to \wb{U_\alpha}$ is closed
   (see e.g. \cite[Thm 3.1.16]{Engelking}), hence $x\in\pi(C)$.\\
   (b) Immediate since $(s\circ\pi)^{-1}([0,\alpha])$ is contained in $s^{-1}([0,\alpha])\times[0,\alpha]$
       which is Lindel\"of.\\
   (c) Let $C\subset \text{Down}(X,s)$ be closed and non-Lindel\"of.
   Then $\pi(C)$ is closed non-Lindel\"of, let $f\co X\to\LL_{\ge 0}$ witness
   the fact that $\pi(C)$ is not f-narrow in $X$. 
   Then $f\circ\pi\co \text{Down}(X,s)\to\LL_{\ge 0}$ witnesses that $C$ is not f-narrow in 
   $\text{Down}(X,s)$.\\
   (d) Let $x=\pi(y)$ for $y\in D$. 
       By definition, $\pi^{-1}(\{x\}) = \{x\}\times[0,s(x)]$. For each $y\in\pi^{-1}(\{x\})$ there 
       is an open $V_y\subset X$ and an interval $(a_y,b_y)\subset[0,s(x)]$ 
       such that $V_y\times(a_y,b_y)$ intersects at most one point of $D$ and contains $y$.
       Then $[0,s(x)]$ is covered by finitely many intervals $(a_y,b_y)$,
       taking the intersection of the corresponding $V_y$ yields a neighborhood of
       $x$ intersecting $\pi(D)$ in at most finitely many points, hence $\pi(D)$ is discrete.\\
   (e) If $C\subset\pi(D)$ is closed discrete, taking one member in $\pi^{-1}(\{c\})\cap D$ for each $c\in C$
       yields a closed discrete subset of $D$.\\
   (f) One direction is immediate.
       A Type I countably compact space is $\omega$-bounded, and $\omega$-boundedness is productive 
       (Lemma \ref{lemma:omegabounded}).
       It follows that $X\times\LL_{\ge 0}$ is $\omega$-bounded, and this holds as well for any closed subspace.\\
   (g) Immediate by (d) and (e).\\
   (h) Local compactness is productive and $\text{Down}(X,s)$ is a closed subspace of $X\times\LL_{\ge 0}$.
\endproof

We now pile up copies of $\text{Down}(X,s)$ as follows.
We define $\mathsf{D}^{(\omega)}(X,s)$ (sometimes shortened as $\mathsf{D}^{(\omega)}$ if there is no ambiguity) 
as the disjoint union of copies 
$\mathsf{D}_n$ ($n\in\omega$) of $\text{Down}(X,s)$ 
where we identify $\langle x,s(x)\rangle\in \mathsf{D}_n$ with 
$\langle x,0\rangle\in \mathsf{D}_{n+1}$. (Notice that $\langle x,s(x)\rangle\in\Delta(X,s)$.
If $X=\LL_{\ge 0}$ and $s$ is the identity, this consists of piling up $\omega$ copies of $\mathbb{O}$,
identifying $\Delta$ in one copy with $\HH$ in the next.)
By a slight abuse of notation we still call $\mathsf{D}_n$ the images of these subspaces
in $\mathsf{D}^{(\omega)}$ (they are thus not disjoint).
The map $s^{(\omega)}\co\mathsf{D}^{(\omega)}\to\LL_{\ge 0}$ defined as $s\circ\pi$ in each
$\mathsf{D}_n$ is a slicer.
We write $\Delta_n$ for the copy of $\Delta(X,s)$ inside $\mathsf{D}_n$.

\vskip .3cm\noindent
We now define $\Omega(X,s)$ to be the disjoint union of $\mathsf{D}^{(\omega)}$
and a copy $\mathsf{L}$ of $\LL_{\ge 0}$ with the following topology. First,
$\mathsf{D}^{(\omega)}$ is open in $\Omega(X,s)$.
Given $a<b$ in $\LL_{\ge 0}$ and $n\in\omega$, set
$$ V_{a,b,n} = \bigl(s^{(\omega)}\bigr)^{-1}\bigl((a,b)\bigr) - \cup_{k<n}\mathsf{D}_k.$$
Then the neighborhoods of $x\in \mathsf{L}$ are $(a,b)\subset\mathsf{L}$ union $V_{a,b,n}$, for $a<x<b$ and $n\in\omega$.
In words: the points in $\mathsf{D}_n$ whose value under $s\circ\pi$ is $x\in\LL_{\ge 0}$
converge to (the copy of) $x\in \mathsf{L}$ when $n$ grows. These points may be thought as being under $x$.
If $A\subset \mathsf{L}$, we write $\Omega(X,A,s)$ for the subspace 
$\mathsf{D}^{(\omega)}\sqcup A$ of $\Omega(X,s)$.
Notice that the function $s_\Omega$ 
defined by the identity on $\mathsf{L}$ and $s^{(\omega)}$ on $\mathsf{D}^{(\omega)}$
is a slicer for $\Omega(X,s)$
(the preimage of $[0,\alpha]$ is a countable union of Lindel\"of spaces
and is thus Lindel\"of), which is thus f-Type I.

\begin{lemma}\label{lemma:Omegaprop}
   Let $X$ be f-Type I with slicer $s\co X\to\LL_{\ge 0}$.
   Then the following hold.\\
   (a) $\Omega(X,s)$ and $\mathsf{D}^{(\omega)}(X,s)$ are Tychonoff iff $X$ is regular.\\
   (b) $\Omega(X,s)$ is countably compact iff $X$ is countably compact.\\
   (c) $\mathsf{D}^{(\omega)}(X,s)$ is $\omega_1$-compact iff $X$ is $\omega_1$-compact.\\
   (d) $\mathsf{D}^{(\omega)}(X,s)$ is locally compact iff $X$ is locally compact.
\end{lemma}
\proof\
   \\
   (a)
   If $X$ is regular, it is straightforward to check that $s_\Omega^{-1}([0,\alpha])$
   is regular for each $\alpha$. Then apply Lemma \ref{lemma:TypeIregular}.\\
   (b) By Lemma \ref{lemma:down} (f) and construction: a countable set either has infinite intersection with some $\mathsf{D}_k$
       or has an accumulation point in $\mathsf{L}$.\\
   (c) An uncountable subset of $\mathsf{D}^{(\omega)}(X,s)$ has an uncountable intersection
       with some $\mathsf{D}_k$, which is a closed subset of $\mathsf{D}^{(\omega)}(X,s)$. 
       The result follows by Lemma \ref{lemma:down} (g).
       \\
   (d) Immediate by Lemma \ref{lemma:down} (h) and construction.
\endproof

Note however that $\Omega(X,A,s)$ is not locally compact at points of $A\subset\mathsf{L}$ 
if the bones of $X$ do not have compact neighborhoods. 
The following lemma is reminiscent of Theorem \ref{thm:stationary} and Lemma \ref{lemma:7.2}.
\begin{lemma}\label{lemma:Anarrow}
   Let $X$ be f-Type I with slicer $s$, and let $B_\alpha = \wb{s^{-1}([0,\alpha))}-s^{-1}([0,\alpha))$
   for $\alpha\in\omega_1$ be the bones of the canonical cover given by $s$.
   If $B_\alpha\not=\varnothing$ for a stationary set of $\alpha\in\omega_1$
   and $A\subset\mathsf{L}$. Then
   $A$ is narrow in $\Omega(X,A,s)$.
\end{lemma}
\proof
   We assume $A$ to be unbounded in $\mathsf{L}$, otherwise the result is immediate.
   For $x\in\mathsf{L}$ and $n\in\omega$, 
   let $N(x,n)=\bigl(s^{(\omega)}\bigr)^{-1}\bigl(\{x\}\bigr) - \cup_{k<n}\mathsf{D}_k$,
   in words: the points of $\mathsf{D}^{(\omega)}$ which are under $x$
   but not in the $n$ first $\mathsf{D}_k$.
   Let $\mathcal{U}=\{U_\alpha\,:\,\alpha\in\omega_1\}$ be a systematic cover of $\Omega(X,A,s)$.
   Suppose that there is some $\alpha$ such that $A_0= \wb{U_\alpha}\cap A$ is non-Lindel\"of.
   If $A-\wb{U_{\alpha+2}}$ is bounded in $\mathsf{L}$ (hence Lindel\"of), 
   then $A$ is contained in $U_\beta$ for some $\beta>\alpha$ and we are over.
   Otherwise, $A_1 = A-U_{\alpha+2}$ is unbounded in $\mathsf{L}$.
   Notice that $A_0\cap A_1=\varnothing$ and both are closed.
   For $i=0,1$,
   $A_i$
   contains an uncountable subset $B_i=\{b^i_\alpha\,:\,\alpha\in\omega_1\}$ such that $b^i_\alpha > \alpha$
   (recall that $A\subset \mathsf{L}$ which is a copy of $\LL_{\ge 0}$). 
   Up to taking subsets, we may assume that there is $n\in\omega$ such that $U_\alpha$ contains 
   $N(b^0_\alpha,n)$ and $\Omega(X,A,s)-\wb{U_{\alpha+1}}$ contains $N(b^1_\alpha,n)$ for each $\alpha$. 
   The closure of $B_i$ in $\mathsf{L}$
   contains a club $C_i\subset\omega_1\subset\mathsf{L}$.
   Let now $C=C_0\cap C_1$.
   Since $C$ is club and $B_\alpha\not=\varnothing$ for a stationary set $S$ of $\alpha$,
   by construction and closedness
   both $\wb{U_\alpha}$ and $\Omega(X,A,s)-U_{\alpha+1}$ contain all of $N(\alpha,n)$ for those $\alpha\in S\cap C$,
   which is clearly impossible.
   It follows that $\wb{U_\alpha}\cap A$
   is Lindel\"of for each $\alpha$, hence $A$ is narrow in $\Omega(X,A,s)$.
\endproof

\begin{lemma}\label{lemma:Omeganonnarrow}
   Let $X$ be f-Type I with slicer $s$ and $f\co\text{Down}(X,s)\to\LL_{\ge 0}$ be 
   such that $f(x)\le s\circ\pi(x)$ for all $x\in X$.
   For each $n\in\omega$ there is $g\co\Omega(X,s)\to\LL_{\ge 0}$
   which is equal to $f$ on $\mathsf{D}_n$ and equal to $s_\Omega$
   on $\mathsf{L}$ and $\mathsf{D}_\ell$ for each $\ell>n+1$.
\end{lemma}
\proof
   Set $h\co\text{Down}(X,s)\to\LL_{\ge 0}$ to be 
   the function $h(\langle x,y\rangle) = \max\{f(\langle x,y\rangle),y\}$.
   Then $h$ is equal to $f$ on $X\times\{0\}$ and to $s\circ\pi$ on $\Delta(X,s)$.
   Define $j\co\text{Down}(X,s)\to\LL_{\ge 0}$ as $j(\langle x,y\rangle) = f(\langle x,0\rangle)$.
   Finally, define $g$ as being equal to $h$ on $\mathsf{D}_{n+1}$, to $f$ in $\mathsf{D}_{n}$,
   to $j$ in $\mathsf{D}_{k}$ with $k<n$ and to $s_\Omega$
   on $\mathsf{L}$ and $\mathsf{D}_\ell$ for each $\ell>n+1$. 
\endproof

\begin{cor}\label{cor:Omeganotfnarrow}
   Let $X$ be f-Type I with slicer $s$ and $A\subset\LL_{\ge 0}$ be unbounded. 
   If $X$ does not contain any non-Lindel\"of closed subspace
   f-narrow in $X$, then any non-Lindel\"of closed subset f-narrow in $\Omega(X,A,s)$ 
   is contained in $A$ outside of $s_\Omega^{-1}([0,\alpha])$ (for some $\alpha$). 
\end{cor}
\proof
   By Lemma \ref{lemma:down} (c), $\text{Down}(X,s)$ has no closed non-Lindel\"of subset f-narrow in it.
   Let $B\subset\Omega(X,A,s)$ be closed.
   If $B\cap\mathsf{D}_k$ is non-Lindel\"of for some $k$, by Theorem \ref{thm:fnarrow}
   there is some $f\co\mathsf{D}_k\to\LL_{\ge 0}$ which is unbounded on $B\cap\mathsf{D}_k$ 
   with $f^{-1}(\{0\})\cap B$ non-Lindel\"of.
   A classical leapfrog argument shows that
   the set of $\alpha$ such that $f((s\circ\pi)^{-1}([0,\alpha]))\subset [0,\alpha]$ is club in $\omega_1$.
   It follows that $h(x) = \min\{s\circ\pi(x),f(x)\}$ is (continuous and) unbounded on $B\cap\mathsf{D}_k$.
   Replacing $f$ by $h$ if needed, we may assume that $f(x)\le s\circ\pi(x)$.
   The function $g$ given by Lemma \ref{lemma:Omeganonnarrow} 
   is unbounded on $B$ and $g^{-1}(\{0\})$ is non-Lindel\"of. Hence $B$ is not f-narrow in $\Omega(X,A,s)$.
   It follows that any closed subspace f-narrow in $\Omega(X,A,s)$ has a Lindel\"of intersection
   with $\mathsf{D}^{(\omega)}$.
\endproof

\begin{example}\label{ex:discreteplus}
   \ \\
   There are f-Type I locally metrizable spaces $Y$ 
   which contain a non-Lindel\"of closed subspace narrow in $Y$
   such that any closed subspace f-narrow in $Y$ is discrete outside of a Lindel\"of subset, and
   with the additional following properties.\\
   (a) Each non-Lindel\"of closed subset of $Y$ contains an uncountable closed discrete subset.
       \\
   (b) Under $\diamondsuit^*$, $Y$ can be made as in (a) and cwH.
   \\
   (c) If there is a Suslin tree, $Y$ can be made to have an
       open dense $\omega_1$-compact and locally compact subspace.
\end{example}
Recall that a tree is {\em Aronszajn} if it is an $\omega_1$-tree whose chains are at most countable.
\proof[Details]
    In each case, $Y=\Omega(T,A,s)$ where $T$ is an Aronszajn tree, $s$ the height function and 
    $A$ the subset of successor ordinals in $\omega_1\subset\LL_{\ge 0}$.
    Then $\Omega(T,A,s)$ is f-Type I and locally metrizable since
    $s_\Omega^{-1}([0,\alpha))$ is second countable for each $\alpha$.
    By construction, $A$ is closed discrete and Lemma \ref{lemma:Anarrow} implies its narrowness in $\Omega(T,A,s)$.
    By Corollary \ref{cor:Omeganotfnarrow} and the fact that an Aronszajn tree does not contain
    a non-Lindel\"of closed f-narrow space (Lemma \ref{lemma:treenotnarrow}),
    any subspace f-narrow in $\Omega(T,A,s)$ is contained in $A$ outside of a Lindel\"of subset.
    This proves the general claim.\\
    (a) Take $T$ to be $\R$-special (that is, there is a non-necesarily continuous
        order preserving map $T\to\R$). $\R$-special Aronszajn trees exist in {\bf ZFC}.
        Any non-Lindel\"of closed subset of an Aronszajn $\R$-special tree contains a 
        uncountable closed discrete subset (see e.g. \cite[Thm 2]{DevlinNoteBaumgartner}).
        Let $C\subset\Omega(T,A,s)$ be closed and non-Lindel\"of.
        If $C\cap\mathsf{D}_k$ is non-Lindel\"of
        for some $k$, by Lemma \ref{lemma:down} (e), $C$ contains an uncountable closed discrete subset.
        If not, its intersection with $A$ is non-Lindel\"of, and we conclude the same.\\
    (b) In \cite{DevlinShelah}, Devlin and Shelah use $\diamondsuit^*$ to construct an $\R$-special
        cwH tree. 
        \begin{claim}
        If $T$ is a cwH $\omega_1$-tree and $A$ is the successor ordinals, then $\Omega(T,A,s)$ is cwH.
        \end{claim}
        \proof
        Let $B$ be closed discrete in $\Omega(T,A,s)$.
        Then $B\cap\mathsf{D}_n$ is closed discrete for each $n$,
        and by Lemma \ref{lemma:down} (d) the projection of $B\cap\mathsf{D}_n$
        on the first factor $T$ of $\mathsf{D}_n$ is closed discrete.
        Since $T$ is cwH (and $s^{-1}([0,\alpha])$ is countable hence ccc), by Lemma \ref{lemma:cwHdiscrete},
        $$
            B\cap\mathsf{D}_n \cap \Bigl(\wb{(s\circ\pi)^{-1}([0,\alpha))} - (s\circ\pi)^{-1}([0,\alpha))\Bigr)=\varnothing
        $$ 
        for a club set $E_n\subset\omega_1$
        of limit ordinals $\alpha$.
        Hence, $B$ avoids $\wb{s_\Omega^{-1}([0,\alpha))}-s_\Omega^{-1}([0,\alpha))$ for $\alpha$ in $E=\cap_{n\in\omega}E_n$.
        Enumerating $E$ as $\{\alpha_\beta\,:\,\beta\in\omega_1\}$,
        $B$ is contained in the disjoint union of the open second countable hence metrizable subspaces 
        $s_\Omega^{-1}\bigl((\alpha_\beta,\alpha_{\beta+1})\bigr)$.
        We may then separate $B$ by open sets inside each piece, yielding the result.
        \endproof
    \noindent
    (c) Take $T$ to be a Suslin tree. Then $\mathsf{D}^{(\omega)}(T,s)$
        is open dense in $\Omega(T,A,s)$, $\omega_1$-compact and locally compact by Lemma \ref{lemma:Omegaprop} (c)--(d).
\endproof

A way to describe (b) is to say that given a subspace narrow in $Y$, there is a skeleton which is totally blind
to it as their intersection is empty, as shown in the proof.\footnote{We remark that 
skeletons being blind is in spectacular accordance with cutting edge research in the
biology of vertebrates.}

%%%%%%%%%%%%%%%%%%%%%%%%%%%%%%%%%%%%%%%%%%%%%%%%%%%%%%%%%%%%%%%%%%%%%%%%%%%%%%%%%%%%%%%%%%%%%%%%%%%%%%%%%%%%%%%%%%%%%%%%%%%%%%%%%%%

\section{Two partial orders on closed (non-Lindel\"of) subsets of a space}\label{sec:orders}

\begin{defi}
   Let $Y$ be a space and $C,D$ be closed subsets of $Y$.
   \\
   (a) Set $C\preceq D$ iff given any systematic cover $\{U_\alpha\,:\,\alpha\in\omega_1\}$ of $Y$,
   if $D\subset U_\alpha$ for some $\alpha$, then $C\subset U_\beta$ for some $\beta$.
   If $C\preceq D\preceq C$, we write $D\equiv C$ and say that $D$ and $C$ are sc-equivalent (sc is for systematic cover).
   \\
   (b) Set
   $C\preceq_f D$ iff for any $f\co Y\to \LL_{\ge 0}$, 
   $f\upharpoonright D$ is bounded implies $f\upharpoonright C$ is bounded.
   If $C\preceq_f D\preceq_f C$, we write $D\equiv_f C$ and say that $D$ and $C$ are uf-equivalent (uf is for unbounded function).
\end{defi}

If $C\preceq D$ and $C\not\equiv D$, we write $C\prec D$, and idem for $\prec_f$.
The order $\preceq_f$ is written $\preceq$ in \cite{mesziguesDirections} and called the UFO,
but being some years older than when \cite{mesziguesDirections} was written enables us to 
slither away from names that sound too cool.
Notice that any closed subset $D$ is sc- and uf-equivalent to the closure of $D-E$, where $E$ is Lindel\"of.
The following lemmas are immediate.

\begin{lemma}
   Let $Y$ be a space and $D\subset Y$ be closed and non-Lindel\"of.
   If $D$ is Type I [resp. f-Type I] in Y, then $D$ is narrow [resp. f-narrow] in $Y$ iff 
      $D\equiv C$ [resp. $D\equiv_f C$] for any closed non-Lindel\"of $C\subset D$.
\end{lemma}

Let $C\subset_L D$ iff $\wb{C-D}$ is Lindel\"of.
\begin{lemma} Let $C,D\subset Y$ be closed.
   Then $$C\subset D\quad\Longrightarrow C\subset_L D\quad\Longrightarrow\quad C\prec D\quad\Longrightarrow\quad C\prec_f D.$$
   Moreover, if $C$ and $D$ are Lindel\"of, then $C\equiv_f D\equiv C$.
\end{lemma}
In particular, the whole space is always $\preceq$ and $\preceq_f$ maximal.

\begin{example}\label{ex:octantorder}\ \\
   (a) In the octant $\mathbb{O}$, all horizontal lines are sc-equivalent and $\prec$ the diagonal, and any non-Lindel\"of
       closed set narrow in it is sc-equivalent to the horizontals or to the diagonal.\\
   (b) In the space of Example \ref{ex:narnotfunc}, $\mathsf{W}\succ\VV_\alpha$ but 
       $\mathsf{W}\equiv_f\VV_\alpha$ when $\alpha<\omega_1$.
\end{example}
\proof[Details]
   (a) follows essentially by Lemma \ref{lemma:propoct} and Example \ref{ex:LxR}, and (b) by construction.
\endproof

Since the whole space is always sc- and uf-maximal, a more interesting poset arises when one looks 
only at (f-)narrow subsets.
\begin{defi}\ \\
   $\mathfrak{N}(D,X)$ is the poset ordered by $\preceq$ of sc-equivalence classes 
   of subsets of $D\subset X$ which are narrow in $X$.\\
   $\mathfrak{N}_f(D,X)$ is the poset ordered by $\preceq_f$ of uf-equivalence classes 
   of subsets of $D\subset X$ which are f-narrow in $X$.
\end{defi}
We denote $\mathfrak{N}(X,X)$ by $\mathfrak{N}(X)$, and the same for $\mathfrak{N}_f(X,X)$.

\begin{lemma}\label{lemma:normal2}
   If $X$ is a normal space and $D$ is closed in $X$, then $\mathfrak{N}(D,X)=\mathfrak{N}_f(D,X)$ as posets.
\end{lemma}
\proof 
   By Lemma \ref{lemma:normalfnar}, a closed subset of $D\subset X$ is narrow iff it is f-narrow.
   The slicers provided by Lemma \ref{lemma:Uknormal} enable to pass from systematic covers to functions
   $X\to\LL_{\ge 0}$ with the same properties, hence $C_0\preceq C_1\:\Longleftrightarrow\: C_0\preceq_f C_1$
   for closed subsets.
\endproof

We show in Lemma \ref{lemma:countclosed} that 
$\mathfrak{N}(X)$ and $\mathfrak{N}_f(X)$ are countably closed when $X$ is Type I and countably compact,
this is not the case if $X$ is only $\omega_1$-compact.

\begin{example}\label{ex:nominmax}
   An $\omega_1$-compact surface (i.e. $2$-manifold) $X$ such that $\mathfrak{N}(X)=\mathfrak{N}_f(X)\simeq\Z$,
   hence $\mathfrak{N}(X)$ has neither a maximal nor a minimal element.
\end{example}
\proof[Details]
    We take copies $O_n$ of $\mathbb{O}$ for $n\in\Z$,
    gluing the copy of $\HH$ in $O_{n+1}$ to that of $\Delta$ in $O_n$.
    More precisely: we identify $\langle x,x\rangle\in O_n$
    with $\langle x,0\rangle\in O_{n+1}$ for each $n\in\Z$.
    (This is very similar to $\mathsf{D}^{(\omega)}(\mathbb{O},\pi)$ where $\pi$ is the projection on the first factor,
    but we take copies ordered by $\Z$ instead of $\omega$.) 
    The resulting space is $\omega_1$-compact since it is a countable union of countably compact subspaces.
    It is not difficult to show that it is also normal, hence 
    $\mathfrak{N}(X)$ and $\mathfrak{N}_f(X)$ agree by Lemma \ref{lemma:normal2}.
    Any closed non-Lindel\"of subset must intersect one of the $O_n$ on a non-Lindel\"of subset, 
    and proceeding as in Example \ref{ex:octantorder} and Lemma \ref{lemma:Omeganonnarrow},
    it is easy to show that the copy of $\Delta$ in $O_n$ is $\prec_f$ the one in $O_{n+1}$,
    and that any non-Lindel\"of narrow subspace is sc-equivalent to one of those diagonals.
\endproof

By adding a copy of $\LL_{\ge 0}$ at the top and/or the bottom, we may continue piling up copies of the octant
and turn $\mathfrak{N}(X)$ to be (for instance) any countable ordinal.
By piling up uncountably many copies of $\mathbb{O}_{\ge\alpha}$ for $\alpha\in\omega_1$,
as in the construction of $S^{(\omega_1,\downarrow)}$ in Example \ref{ex:nomin} below,
we may obtain $\omega_1$, or actually any ordinal $<\omega_2$.
(The resulting space 
is countably compact exactly when the ordinal is successor, $\omega_1$-compact when it has countable cofinality,
not $\omega_1$-compact when the ordinal has uncountable cofinality.)
Hence, it is not completely trivial to come up with a countably compact space such that $\mathfrak{N}(X)$ has no maximum or minimum,
and our next results are geared towards this precise problem.
We first state the following easy lemma whose proof is left to the reader (or see \cite[Prop. 2.2]{EisworthNyikos:2005}).
\begin{lemma}\label{lemma:capclnL}
   If $C_0\supset C_1\supset C_2\supset\dots$ is an $\omega$-sequence of closed non-Lindel\"of sets
   of a countably compact Type I space, then $\cap_{n\in\omega}C_n$ is closed and non-Lindel\"of.
\end{lemma}

\begin{lemma}\label{lemma:countclosed}
   If $X$ is Type I and countably compact, then
   $\mathfrak{N}(X)$ and $\mathfrak{N}_f(X)$ are up- and downwards countably closed.
\end{lemma}
\proof
   Let $D_n$ be narrow in $X$ for $n\in\omega$.
   We may assume that $D_n$ is non-Lindel\"of for each $n$. 
   Let 
   $$ D = \bigcap_{n\in\omega}\wb{\bigcup_{k\ge n} D_k}.$$
   Then $D$ is closed and non-Lindel\"of by Lemma \ref{lemma:capclnL}.
   \begin{claim}\label{claim:narlim}
      If $\{U_\alpha\,:\,\alpha\in\omega_1\}$ is a systematic cover of $X$ such that
      $\wb{U_0}\cap D$ is non-Lindel\"of, then there is strictly increasing sequence 
      $k_n$ ($n\in\omega$) and $\alpha\in\omega_1$
      such that $U_\alpha\supset D_{k_n}$ for each $n$.
   \end{claim}
   \proof
      By induction on $n$. Let $k_{-1}=\alpha_{-1}=0$. 
      Since the closure of Lindel\"of subspaces
      are Lindel\"of in Type I spaces, 
      if $\wb{U_{\alpha_n+1}}\cap D_{m}$ is Lindel\"of for each $m>k_n$, 
      then $D\cap\wb{U_0}\subset\wb{U_{\alpha_n+1}}\cap\wb{\bigcup_{m>k_n}D_m}$ would be Lindel\"of.
      Hence there is $k_{n+1}\ge k_n$ such that $\wb{U_{\alpha_n+1}}\cap D_{k_{n+1}}$ is non-Lindel\"of,
      and then $D_{k_{n+1}}\subset U_{\alpha_{n+1}}$ for some $\alpha_{n+1}$ since $D_{k_{n+1}}$ is narrow in $X$.
      Set $\alpha$ to be the supremum of the $\alpha_n$ to conclude.
   \endproof
   \noindent Assume now that $D_n\prec D_{n+1}$ for each $n$. Then 
   $D\succeq D_n$ for each $n\ge 0$ by the previous claim and transitivity of $\preceq$.
   Actually, the claim also shows that $D$ is narrow in $X$, because if the intersection
   of $D$ with the closure of a member of a systematic cover $\mathcal{U}$ is non-Lindel\"of,
   then each $D_n$ lies in some member of $\mathcal{U}$, and $D$ is then inside the closure of
   their union.
   To finish, it is enough to show that $D\not\preceq D_n$ for each $n$.
   For this, let $\{U_\alpha\,:\,\alpha\in\omega_1\}$ be a systematic cover of $X$
   witnessing that $D_n\prec D_{n+1}$, that is, $D_n\subset U_0$ but 
   $D_{n+1}\cap\wb{U_\alpha}$ is Lindel\"of for each $\alpha$.
   It follows that $D_{k}\cap\wb{U_\alpha}$ is Lindel\"of for each $k\ge n+1$
   by transitivity of $\preceq$.
   Again, the closure of a Lindel\"of set is Lindel\"of, hence
   $D\cap \wb{U_\alpha}\subset \wb{\cup_{k\ge n+1} D_k}\cap\wb{U_\alpha}$ is Lindel\"of as well.
   This shows that $D$ is an $\preceq$-upper bound of the $D_n$.\\
   Now, if $D_n\succ D_{n+1}$ for each $n$, fix $n$ and let $\mathcal{U}=\{U_\alpha\,:\,\alpha\in\omega_1\}$ be a systematic cover 
   such that $U_0\supset D_{n+1}$. 
   Then, there is $\alpha$ such that $D_k\subset U_\alpha$ for each $k>n$.
   Since $D\subset \wb{\cup_{k\ge n+1} D_k}\cap\wb{U_\alpha}$, $D\preceq D_n$.
   If $\mathcal{U}$ is such that moreover
   $\wb{U_\alpha}\cap D_{n}$ is Lindel\"of for each $\alpha$, then we see that $D\not\succeq D_{n}$, and it follows that 
   $D\prec D_n$ for each $n$. The previous claim easily shows that
   $D$ is narrow in $X$, which concludes the proof for $\mathfrak{N}(X)$.
   The proof for $\mathfrak{N}_f(X)$ is the same, replacing systematic covers with 
   appropriate functions $X\to\LL_{\ge 0}$.
\endproof

It is however not very difficult to find an
example without a $\preceq$-minimal element.

\begin{figure}
  \begin{center}
    \epsfig{figure = 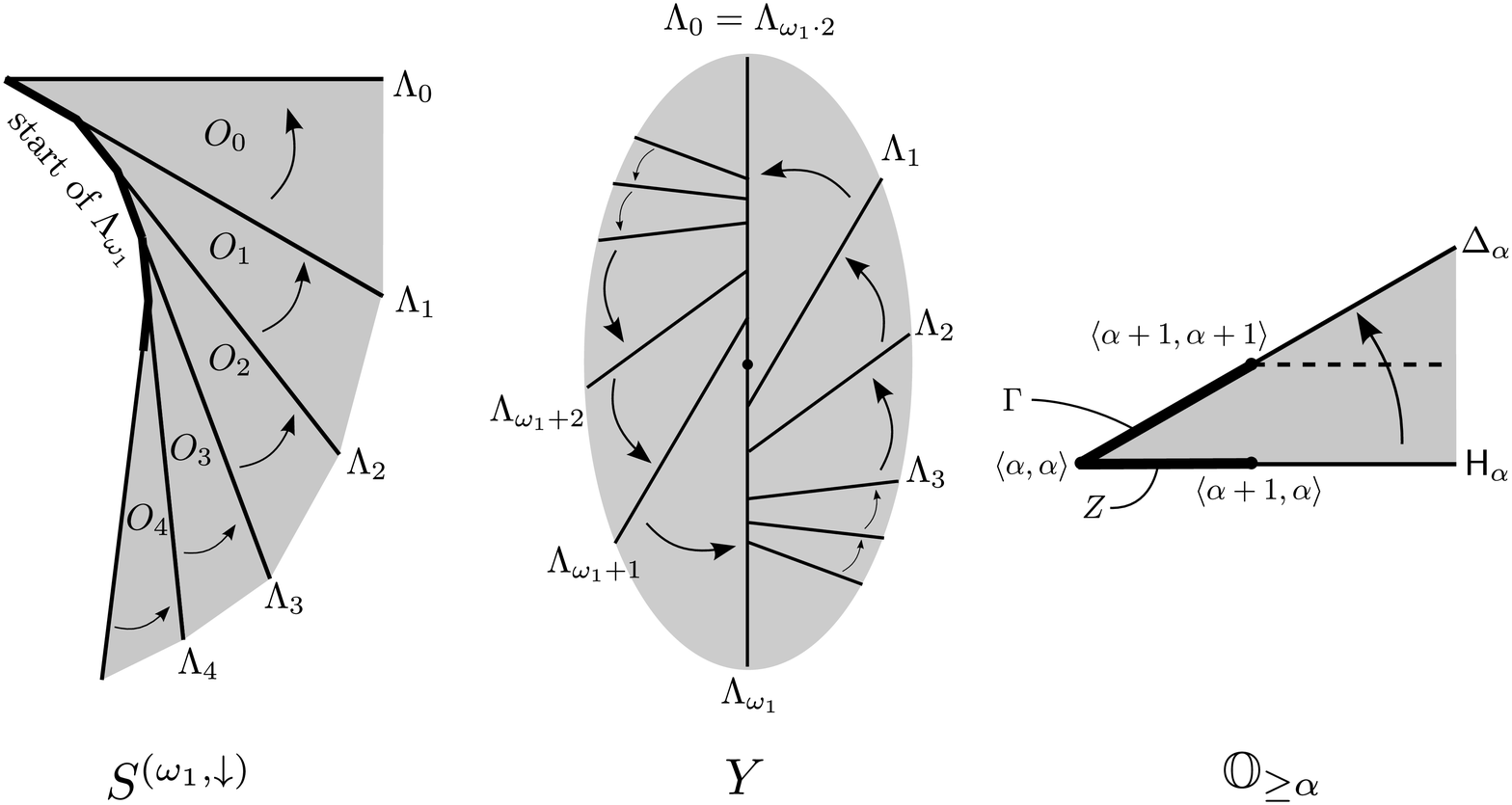, width = .85\textwidth}
    \caption{Example \ref{ex:nomin}.}
    \label{fig:nomin}
  \end{center}
\end{figure}

\begin{example}\label{ex:nomin}
   A Type I countably compact surface $Y$ such that $\mathfrak{N}(Y) = \mathfrak{N}_f(Y)$ has no minimal element.
\end{example}
\proof[Details]
    The construction is very similar to that of Examples \ref{ex:narnotfunc} and \ref{ex:nominmax}
    and is shown in Figure \ref{fig:nomin} (center).
    We start with piling up (non disjoint) copies $O_\alpha$ of $\mathbb{O}_{\ge\alpha}$ for $\alpha\in\omega_1$, 
    gluing $\Delta_{\alpha+1}$ to $\HH_\alpha$, that is,
    we identify $\langle x,x\rangle\in\Delta_{\alpha+1}$ with $\langle x,0\rangle\in \HH_{\alpha}$
    (in a sense, we put $\mathbb{O}_{\ge\alpha+1}$ below $\mathbb{O}_{\ge\alpha}$).
    At limit levels we define neighborhoods similarly as in Example \ref{ex:narnotfunc}:
    points with first coordinate $x\in\LL_{\ge 0}$ in $O_\alpha$ for $\alpha<\beta$
    with limit $\beta$ accumulate to $\langle x,x \rangle \in O_\beta$. 
    When $\alpha$ is limit or $0$,
    we let $\Lambda_\alpha$ be the copy of $\Delta_\alpha$ in $O_\alpha$.
    For successors, we define $\Lambda_{\alpha+1}$ to be the copy of $\HH_\alpha$ in $O_\alpha$
    (which is identified with
    $\Delta_{\alpha+1}\subset O_{\alpha+1}$ above $\alpha+1$).
    Then, as in Example \ref{ex:nominmax}, $\Lambda_\alpha \succ \Lambda_\beta$ if
    $\alpha<\beta$. Call $S^{(\omega_1,\downarrow)}$ the space thus obtained,
    whose first stages of construction are depiced in Figure \ref{fig:nomin} (left).
    Set $\Lambda_{\omega_1}$ to be the union of the segments 
    $\{\langle y,y\rangle\in\Lambda_\alpha\,:\,\alpha\le y\le \alpha+1\}$.
    Then $\Lambda_{\omega_1}$ is a copy of $\LL_{\ge 0}$ 
    uf-uncomparable (and thus, sc-uncomparable) with each $\Lambda_\alpha$.
    To see this, let $f_\alpha\co\mathbb{O}_{\ge\alpha}\to[0,\alpha+1]\subset\LL_{\ge 0}$ be a continuous map which
    takes value $\alpha$ on the segment $\Gamma$ and above the dashed line depicted in Figure \ref{fig:nomin} (right),
    is equal to the projection on the first factor on the segment $Z$
    and constant with value $\alpha+1$ on the bottom boundary $\HH_\alpha$ above $\langle\alpha+1,0\rangle$.
    Then $f\co S^{(\omega_1,\downarrow)}\to\LL_{\ge 0}$
    defined as $f\upharpoonright O_\alpha = f_\alpha$ is well defined, continuous, bounded on each $\Lambda_\alpha$
    when $\alpha<\omega_1$ and unbounded on $\Lambda_{\omega_1}$, showing that $\Lambda_\alpha\not\succeq_f\Lambda_{\omega_1}$.
    Conversely, the map $g$ defined as constant on $\alpha$ on $O_\beta$ when $\alpha<\beta$,
    to the projection on the second coordinate on $O_\alpha$ and to the projection on the first coordinate
    on $O_\beta$ with $\alpha>\beta$ is continuous, unbounded on $\Lambda_\alpha$ and bounded on $\Lambda_{\omega_1}$.
    It follows that $\Lambda_\alpha\not\preceq_f\Lambda_{\omega_1}$. 
    Notice that any closed non-Lindel\"of subset of $S^{(\omega_1,\downarrow)}$ intersects
    either one of the $O_\alpha$ or $\Lambda_{\omega_1}$ in a non-Lindel\"of set
    (this is similar to Claims \ref{claim:3.10.1}--\ref{claim:3.10.2}).
    \\
    Now, define $Y$ by gluing two copies of $S^{(\omega_1,\downarrow)}$ by 
    identifying pointwise $\Lambda_0$ in the first copy with $\Lambda_{\omega_1}$
    in the second, and vice-versa. To simplify, shift the indexing in the second copy by $\omega_1$, as in Figure     
    \ref{fig:nomin} (middle). Then any non-Lindel\"of space narrow in $Y$ must be sc-equivalent to $\Lambda_\alpha$
    for some $\alpha < \omega_1\cdot 2$.
    It is easy to adapt the maps $S^{(\omega_1,\downarrow)}\to\LL_{\ge 0}$ witnessing that $\Lambda_{\omega_1}$ and $\Lambda_\alpha$
    are uncomparable to have domain all of $Y$.
    It follows that $\mathfrak{N}_f(Y)$ is the partially ordered set consisting of two disjoint
    copies of $\omega_1$ with reverse order, elements in distinct copies being uncomparable.
    One may check that $Y$ is a normal space, 
    hence $\mathfrak{N}(Y)=\mathfrak{N}_f(Y)$, and we are over.
\endproof

On the other hand, we have:
\begin{thm}[{\bf PFA}]
   \label{thm:PFAmax}
   In a Type I countably tight countably compact space $X$, $\mathfrak{N}(X)$ and $\mathfrak{N}_f(X)$ are $\omega_1$-upwards closed.
\end{thm}
\proof
   As Lemma \ref{lemma:countclosed}, we prove it for 
   $\mathfrak{N}(X)$, the proof for $\mathfrak{N}_f(X)$ being basically the same.
   {\bf PFA} is needed only through Theorem \ref{thm:PFA1} at the very end of the proof.
   Let $\mathcal{D}=\{D_\alpha\subset X\,:\,\alpha\in\omega_1\}$ be a familiy 
   of closed subspaces narrow in $X$, with $D_\alpha\prec D_\beta$ whenever $\alpha<\beta$.
   Let $\mathcal{U}=\{U_\alpha\,:\,\alpha\in\omega_1\}$ be a canonical cover for $X$.
   By Lemma \ref{lemma:skeleton}, each $D_\beta$ intersects the members of $\BO(\mathcal{U})$ on a club set $C_\beta$ of indices.
   Set $C$ to be the diagonal intersection of the $C_\beta$. 
   This means that if $\alpha\in C$, then $D_\beta\cap (\wb{U_\alpha}-U_\alpha)\not=\varnothing$
   for each $\beta<\alpha$. We may assume that $C$ contains only limit ordinals.
   We then set:
   \begin{equation}\label{eq:defD}
   \bigtriangleup\mathcal{D} = \wb{\bigcup_{\alpha\in C} H_\alpha},\text{ where } 
   H_\alpha = (\wb{U_\alpha}-U_\alpha)\cap\bigcap_{\beta<\alpha}\wb{\bigcup_{\gamma:\beta<\gamma<\alpha} D_\gamma} .
   \end{equation}
   Then $\bigtriangleup\mathcal{D}$ is a closed unbounded subset 
   of the skeleton and is thus non-Lindel\"of, and each $H_\alpha$ is compact.
   Notice that setting $\mathcal{D}_\delta = \{D_\beta-U_\delta\,:\,\delta<\beta<\omega_1\}$ for $\delta\in\omega_1$,
   have $\bigtriangleup\mathcal{D}_\delta = \bigtriangleup\mathcal{D}-U_\delta$.
   \begin{claim}\label{claim:om1bd1}
      Let $\{V_\alpha\,:\,\alpha\in\omega_1\}$ be a systematic cover of $X$.
      If $\wb{V_\beta}\cap \bigtriangleup\mathcal{D}$ is non-Lindel\"of for some $\beta$, then for each
      $\alpha\in\omega_1$, there is $\beta(\alpha)$ such that $\wb{V_{\beta(\alpha)}}\supset D_\alpha$.
   \end{claim} 
   \proof
      By construction, $\wb{V_{\beta+1}}$ intersects $H_\alpha$ for $\alpha$ in a club $B\subset C$.
      It follows that for each $\alpha\in B$,
      $\wb{V_{\beta+2}}$ intersects $(\wb{U_\alpha}-U_\alpha)\cap D_{\gamma(\alpha)}$ for some 
      $\gamma(\alpha)<\alpha$. Fodor's Lemma then implies that there is $\gamma\in\omega_1$ such that
      $\wb{V_{\beta+2}}$ intersects $(\wb{U_\alpha}-U_\alpha)\cap D_{\gamma}$ for (stationary many, hence) club
      many $\alpha$. Since $D_\gamma$ is narrow in $X$, there is some $\delta$ such that $\wb{V_\delta}\supset D_\gamma$,
      and the same holds for each ordinal $\gamma'$ smaller than $\gamma$ since $D_\gamma\succ D_{\gamma'}$.
      By replacing $\mathcal{D}$ with $\mathcal{D}_{\delta'}$, 
      we may ensure that $\gamma>\delta'$ for any $\delta'<\omega_1$,
      proving the claim.
   \endproof
   By definition of $\preceq$,
   it follows that if 
   $E\subset \bigtriangleup\mathcal{D}$ is closed and non-Lindel\"of, then 
   $E\succeq D_\alpha$ for each $\alpha\in\omega_1$.
   Moreover, $E\not\preceq D_\alpha$ for each $\alpha\in\omega_1$.
   Indeed, since $D_\alpha\prec D_{\alpha+1}$, 
   there is a systematic cover $\{V_\alpha\,:\,\alpha\in\omega_1\}$ of $X$
   with $D_\alpha\subset V_0$ and $D_{\alpha+1}\cap\wb{V_\gamma}$ Lindel\"of for each $\gamma$.
   The previous claim implies that $ \bigtriangleup\mathcal{D}\cap\wb{V_\gamma}$ is Lindel\"of for each $\gamma$,
   and thus so is $\wb{V_\gamma}\cap E$. \\
   Since $\bigtriangleup\mathcal{D}$ is closed non-Lindel\"of in $X$, by Theorem \ref{thm:PFA1}
   there is a copy of $\omega_1$ in it. This copy is an $\preceq$-upper bound of the $D_\alpha$
   in $\mathfrak{N}(X)$.
\endproof

\begin{rem}
   We use {\bf PFA} only to ensure that $ \bigtriangleup\mathcal{D}$ contains a 
   closed non-Lindel\"of subset narrow in itself, and thus in
   $X$ since $X$ is of Type I. 
\end{rem}

We were
frustratingly unable to find the answer to the following question (but have a nagging sensation
of overlooking something simple):
\begin{q}
   Is there a ``reasonable'' (i.e. countably tight, first countable, etc) Type I space $X$
   such that $\mathfrak{N}(X)$ (and/or $\mathfrak{N}_f(X)$) contains a chain of cofinality $\omega_2$ (that is, a
   collection $\{D_\alpha\,:\,\alpha\in\omega_2\}$
   with $D_\alpha\prec D_\beta$ when $\alpha<\beta$)~?
\end{q}
(Actually, we do not know the answer even for unreasonable spaces.) A negative answer to this question
implies that {\bf PFA} ensures (through Zorn's Lemma, Lemma \ref{lemma:countclosed} and Theorem \ref{thm:PFAmax}) 
that $\mathfrak{N}(X)$ has a maximal element when 
$X$ is Type I, countably compact and countably tight.

Let us end this paper with an example showing that $\mathfrak{N}(M)$ can be of cardinality $\mathfrak{c}$.

\begin{figure}
  \begin{center}
    \epsfig{figure = 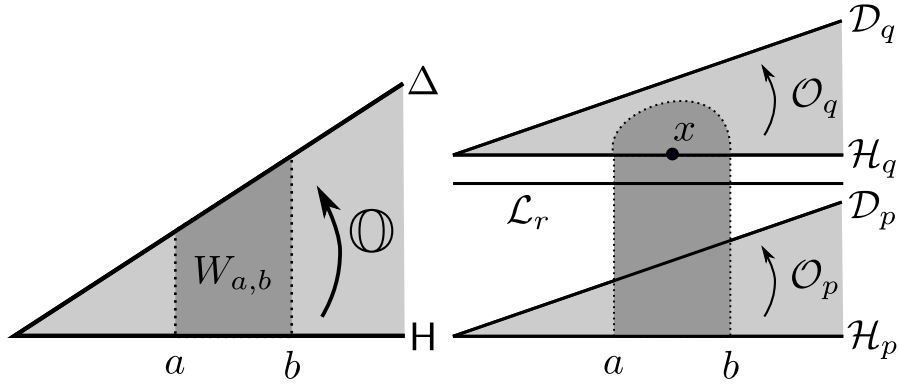, width = .65\textwidth}
    \caption{Example \ref{ex:Cantor}. On the righthandside, $p<r<q$.}
    \label{fig:Cantor}
  \end{center}
\end{figure}

\begin{example}\label{ex:Cantor}
   A Type I countably compact surface $M$ such that $\mathfrak{N}(M)$ is the Cantor set.
\end{example}
\proof[Details]
   Start with $\LL_{\ge 0}\times I$ where $I\subset\R$ is a compact interval with irrational boundary points, delete
   $\LL_{\ge 0}\times\{q\}$ for rational $q$ and replace it by a copy $\mathcal{O}_q = \mathbb{O}\times\{q\}$ of the octant. 
   We thus write points of $M$ as pairs $\langle x,r\rangle$ where $x$ is in $\LL_{\ge 0}$
   or $\mathbb{O}$ according to whether $r$ is irrational or rational.
   Write $\mathcal{D}_q,\mathcal{H}_q$ for 
   the subsets $\Delta\times\{q\},\HH\times\{q\}$, and $\mathcal{L}_r$ for $\LL_{\ge 0}\times\{r\}$
   (with $q$ rational and $r$ irrational, obviously).
   Given $a,b\in\LL_{\ge 0}$ and $r\in[0,1]$, set $W_{a,b}$ to be $\{\langle x,y\rangle\in\mathbb{O}\,:\,a<x<b\}$
   (see Figure \ref{fig:Cantor}, left) and $N_{a,b,r}$ 
   be $(a,b)\times\{r\}\subset \mathcal{L}_r$ if $r$ is irrational and 
   and $W_{a,b}\times\{r\}\subset \mathcal{O}_r$ if $r$ is rational.
   A neighborhood of $\langle x,r\rangle$ for irrational $r$ is given by 
   $\cup_{r'\in(s,t)} N_{a,b,r'}$, for some $a<x<b$ in 
   $\LL_{\ge 0}$ and $s<r<t$ in $[0,1]$ (with the obvious subtleties when $r$ is in the boundary of $I$).
   Fix a rational $q\in I$.
   Neighborhoods of points in the interior of each $\mathcal{O}_q$ (that is: away from
   $\mathcal{D}_q,\mathcal{H}_q$)
   are those of $\mathbb{O}$. Those of points in $\mathcal{H}_q$
   are given by $\cup_{r'\in(s,q)} N_{a,b,r'}$
   for some $s<q$ union $U\times\{q\}$, for $U$ open in $\mathbb{O}$ whose intersection with $\HH$
   is an interval $(a,b)\times\{0\}$ containing $x$. See Figure \ref{fig:Cantor}, right.
   Neighborhoods of points in $\mathcal{D}_q$ are defined similarly, taking $\cup_{r'\in(q,t)} N_{a,b,r'}$
   for some $t>q$.
   \\  
   Let $\tau\co M\to\LL_{\ge 0}\times I$ be the map obtained by sending 
   $\bigl\langle \langle x,y\rangle,q\bigr\rangle$ to 
   $\langle x,q\rangle$ when $q$ is rational and the identity on $\mathcal{L}_r$ when $r$ is irrational.
   Let $\pi\co \LL_{\ge 0}\times I\to\LL_{\ge 0}$ be the projection on the first factor.
   Then $\pi\circ\tau$ is a slicer for $M$, which is easily seen to be a countably compact f-Type I surface (with boundary).
   \begin{claim}
      $\tau$ is a perfect closed quotient map.
   \end{claim}
   \proof The fact that it is perfect (i.e. preimages of points are compact) and a quotient map
          should be clear by definition. Any perfect map with range a locally compact space is closed,
          hence so is $\tau$. 
   \endproof
   Therefore, if $C\subset M$ is closed and non-Lindel\"of, so is $\tau(C)$, and by Lemma \ref{lemma:propoct} (b)
   its intersection with some horizontal line in $\LL_{\ge 0}\times I$ is non-Lindel\"of.
   Hence $C$ has a non-Lindel\"of intersection with either $\mathcal{L}_r$ for some irrational $r$
   or $\mathcal{O}_q$ for some rational $q$, and 
   it follows from the proof of Claim \ref{claim:lastone} below that a closed non-Lindel\"of subspace narrow in $M$
   is sc- and uf-equivalent to $\mathcal{L}_r$, $\mathcal{H}_r$ or $\mathcal{D}_r$ for some $r\in I$.
   \begin{claim}
      Let $U\subset M$ be open, $r$ be irrational and $q$ rational.
      If $U\supset \mathcal{L}_r$, then $\tau(U)\supset \LL_{\ge 0}\times(t,s)$
      for some $t<r<s$.
      If $U\supset \mathcal{H}_q$, then $\tau(U)\supset \LL_{\ge 0}\times(t,q]$ for some $t<q$.
      If $U\supset \mathcal{D}_q$, then $\tau(U)\supset \LL_{\ge 0}\times[q,t)$ for some $t>q$.
   \end{claim}\label{claim:Cantor2}
   \proof
      Suppose that $U\supset \mathcal{L}_r$ and $\tau(U)\not\supset \LL_{\ge 0}\times(t,r]$ for any $t<r$.
      Then there is an $\omega$-sequence $\langle x_n,r_n\rangle\not\in \tau(U)$ with $r_n\nearrow r$.
      Let $x\in\LL_{\ge 0}$ be an accumulation point of the $x_n$.
      By the previous claim
      $\langle x,r\rangle\in \tau(M-U)$, which is impossible since 
      $\tau(M-U)\cap \LL_{\ge 0}\times\{r\}=\varnothing$. The same holds if 
      $\tau(U)\not\supset \LL_{\ge 0}\times[r,s)$ for any $r<s$, with obvious changes if $r$ is a boundary point of $I$.
      \\
      Suppose now that $U\supset \mathcal{H}_q$ and $\tau(U)\not\supset \LL_{\ge 0}\times(t,q]$ for any $t<q$.
      Forget all of $\mathcal{O}_q$ except $\mathcal{H}_q$, and every point with vertical coordinate $>q$
      and argue as above to obtain the desired contradiction. The proof for $\mathcal{D}_q$ is the same.
   \endproof  
   \begin{claim}\label{claim:lastone}
      Let $t<p<r$, with $t,r$ irrational and $p$ rational.
      Then 
      \begin{equation}\label{eq:prec} \mathcal{L}_t \prec \mathcal{H}_p \prec \mathcal{D}_p \prec \mathcal{L}_r. \end{equation}
   \end{claim}
   \proof
   Example \ref{ex:octantorder} shows that $\mathcal{H}_p\preceq\mathcal{D}_p$.
   Let $R_s$ denote $\mathcal{L}_s$ if $s$ is irrational and $\mathcal{H}_s$ if it is rational.
   We show that $R_s\succeq \mathcal{D}_p,\mathcal{H}_q,\mathcal{L}_t$ whenever $t,p<s$, which implies
   (\ref{eq:prec}) with $\preceq$ instead of $\prec$.
   Let $\{U_\alpha\,:\,\alpha\in\omega_1\}$
   be a systematic cover of $M$ such that $U_0\supset R_s$.
   By Claim \ref{claim:Cantor2}, $\tau(U_0)$ contains $\LL_{\ge 0}\times(u,s]$ for some $u<s$.
   If $u$ is rational, it implies by construction that $\wb{U_0}\supset \mathcal{D}_u$.
   Hence, by Example \ref{ex:octantorder}, $\mathcal{H}_u\subset U_\gamma$ for some $\gamma_0$.
   If $u$ is irrational, then $U_1\supset\wb{U_0}\supset \mathcal{L}_u$.
   Proceeding by induction, we obtain a strictly decreasing sequence $u_\alpha\in I$ and ordinals $\gamma_\alpha$
   such that $U_{\gamma_\alpha}$ contains all the points of $M$ with vertical coordinate between $t_\alpha$
   and $s$. The sequence $u_\alpha$ must reach the lower boundary point of $I$ in at most countably many steps,
   hence all $\mathcal{D}_q,\mathcal{H}_p,\mathcal{L}_t$ are contained in some $U_\gamma$ whenever $t,p<s$.
   \\
   To pass from $\preceq$ to $\prec$, 
   notice that it is easy to define a map $f_p\co M\to\LL_{\ge 0}$ which is identically
   $0$ on each $\mathcal{L}_t,\mathcal{O}_t$
   for $t<p$ and on $\mathcal{H}_p$ and equal to $\pi\circ\tau$ on $\mathcal{D}_p$ and $\mathcal{L}_r,\mathcal{O}_r$ for $r>p$
   (take the vertical projection in $\mathcal{O}_p$).
   This shows that $\mathcal{L}_t,\mathcal{H}_p\not\succeq_f \mathcal{D}_p \not\succeq_f \mathcal{L}_r$.
   The remaining $\mathcal{L}_t\not\succeq_f \mathcal{H}_p$ is provided by $f_q$ for some $t<q<p$.
   \endproof
   Hence, $\mathfrak{N}(M)$ is the ordered set obtained by doubling each rational point in $I$,
   with no point in between new pairs.
   This is a way of describing the Cantor space.
\endproof

If we insert copies of $(\LL_{\ge 0})^2$ instead of copies of $\mathbb{O}$
in the previous example, then $\mathfrak{N}(M)$ is the poset obtained by taking $\R$ and replacing 
each rational $q$ by a triplet $q_0,q_1,q_2$, where $q_0,q_2<q_1$, $q_0\perp q_2$,
and each other pair of points are uncomparable. In particular, $\mathfrak{N}(M)$
contains an antichain of cardinality $\mathfrak{c}$ (the irrational points).
Notice also that in Example \ref{ex:Cantor}, countably many octants were inserted at once.
We may instead proceed inductively by removing (final parts of) copies of $\LL_{\ge 0}$ 
above
level $\alpha$ and 
insert copies of $\mathbb{O}_{\ge\alpha}$ in the scar to obtain more complicated posets. 
It is probably possible to 
use Example \ref{ex:Cantor} as a blueprint to obtain some kind of general theorem 
about inverse limits of spaces obtained this way, thus producing many other examples.
We did not pursue this road (except for a small number of specimens described with 
variable rigor in \cite{mesziguesDirections})
and let the inquisitive reader forge their way through this swamp of possibilities
if it suits their curiosity.

\vskip .3cm
{\bf Acknowledgements.}
Some of the work presented here dates back to 2006, when David Gauld invited the author to Auckland's university. We 
express our gratitude to him and this institution for their hospitality. We also thank P. Nyikos for sharing
some of his insights and ideas about the subject, and Jean-Luc le T\'enia for his tenacity.


\begin{thebibliography}{10}

\bibitem{mesziguessurf}
M.~Baillif.
\newblock Homotopy in non-metrizable $\omega$-bounded surfaces.
\newblock preprint arXiv:math/0603515v2, 2006.

\bibitem{mesziguesDirections}
M.~Baillif.
\newblock Directions in {T}ype {I} spaces.
\newblock preprint arXiv:1404.1398, 2014.

\bibitem{Balogh:1989}
Z.~Balogh.
\newblock On compact {H}ausdorff spaces of countable tightness.
\newblock {\em Proc. Amer. Math. Soc.}, 105:755--764, 1989.

\bibitem{DevlinNoteBaumgartner}
K.~J. Devlin.
\newblock Note on a theorem of {J. Baumgartner}.
\newblock {\em Fund. Math}, 76:255--260, 1972.

\bibitem{DevlinShelah}
K.~J. Devlin and S.~Shelah.
\newblock Souslin properties and tree topologies.
\newblock {\em Proc. London Math. Soc.}, 39:237--252, 1979.

\bibitem{Eisworth:2002}
T.~Eisworth.
\newblock On perfect pre-images of $\omega_1$.
\newblock {\em Topology Appl.}, 125(2):263--278, 2002.

\bibitem{EisworthNyikos:2005}
T.~Eisworth and P.~Nyikos.
\newblock First countable, countably compact spaces and the continuum
  hypothesis.
\newblock {\em Trans. Amer. Math. Soc.}, 357:4329--4347, 2005.

\bibitem{EisworthNyikos}
T.~Eisworth and P.~Nyikos.
\newblock Antidiamonds principles and topological applications.
\newblock {\em Trans. Amer. Math. Soc.}, 361:5695--5719, 2009.

\bibitem{Engelking}
R.~Engelking.
\newblock {\em General topology}.
\newblock Heldermann, Berlin, 1989.
\newblock Revised and completed edition.

\bibitem{Gauldbook}
D.~Gauld.
\newblock {\em Non-metrisable manifolds}.
\newblock Springer-Verlag, New York-Berlin, 2014.

\bibitem{GillmanJerisonBook}
L.~Gillman and M.~Jerison.
\newblock {\em Rings of continuous functions}.
\newblock Van Nostrand, Amsterdam, 1960.

\bibitem{Nyikos:1984}
P.~Nyikos.
\newblock The theory of nonmetrizable manifolds.
\newblock In {\em Handbook of {S}et-{T}heoretic {T}opology (Kenneth Kunen and
  Jerry~E. Vaughan, eds.)}, pages 633--684. North-Holland, Amsterdam, 1984.

\bibitem{Nyikos:trees}
P.~Nyikos.
\newblock Various topologies on trees.
\newblock In P.R. Misra and M.~Rajagopalan, editors, {\em Proceedings of the
  Tennessee Topology Conference}, pages 167--198. World Scientific, 1997.

\bibitem{Nyikos:2003}
P.~Nyikos.
\newblock Applications of some strong set-theoretic axioms to locally compact
  {$T_5$} and hereditarily {scwH} spaces.
\newblock {\em Fund. Math.}, 178:25--45, 2003.

\bibitem{CEIT}
L.~A. Steen and J.~A.~Seebach Jr.
\newblock {\em Counterexamples in topology}.
\newblock Springer Verlag, New York, 1978.

\bibitem{vanMill:1984}
J.~van Mill.
\newblock An introduction to $\beta\omega$.
\newblock In {\em Handbook of {S}et-{T}heoretic {T}opology (Kenneth Kunen and
  Jerry~E. Vaughan, eds.)}, pages 503--568. North-Holland, Amsterdam, 1984.

\end{thebibliography}
\end{document}